\documentclass{amsart}
\usepackage{amsfonts}
\usepackage{amssymb}
\usepackage{graphicx}
\usepackage{psfrag}

\newcommand{\N}{{\mathbb N}}
\newcommand{\Z}{{\mathbb Z}}

\newcommand{\C}{{\mathbb C}}

\newcommand{\Groth}{{\mathfrak{G}}}
\newcommand{\bull}{{\scriptscriptstyle \bullet}}
\newcommand{\rank}{\operatorname{rank}}

\newcommand{\Hom}{\operatorname{Hom}}
\newcommand{\Gr}{\operatorname{Gr}}

\newcommand{\OO}{{\mathcal O}}

\newcommand{\Hecke}{{\mathcal H}}
\newcommand{\D}{{\mathfrak D}}
\newcommand{\cross}{\includegraphics[scale=0.1]{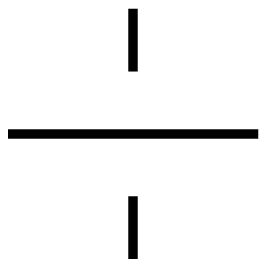}}
\newcommand{\avoid}{\includegraphics[scale=0.1]{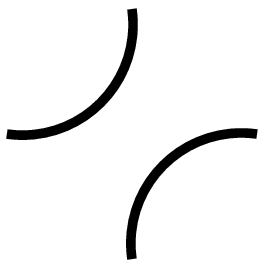}}
\newcommand{\wh}{\widehat}

\newcommand{\tto}{\twoheadrightarrow}

\newcommand{\pic}[2]{\includegraphics[scale=0.#1]{#2.eps}}

\newcommand{\picT}[1]{\includegraphics[scale=0.66]{#1.eps}}

\newtheorem{thm}{Theorem}[section]
\newtheorem{lemma}[thm]{Lemma} 
\newtheorem{prop}[thm]{Proposition}
\newtheorem{cor}[thm]{Corollary}

\theoremstyle{definition}
\newtheorem*{remark}{Remark}

\newcommand{\refthm}[1]{Theorem~\ref{#1}}

\newcommand{\reflemma}[1]{Lemma~\ref{#1}}
\newcommand{\refcor}[1]{Corollary~\ref{#1}}

\newcommand{\extra}[1]{}

\newcommand{\ignore}[1]{}

\psfrag{alpha}{$\alpha$}
\psfrag{Beta}{$\beta$}
\psfrag{w}{$w$}
\psfrag{P1^}{$\wh{P_1}$}
\psfrag{P2^}{$\wh{P_2}$}
\psfrag{Pn^}{$\wh{P_n}$}
\psfrag{x}{$x$}
\psfrag{y}{$y$}
\psfrag{z}{$z$}
\psfrag{a1}{$a_1$}
\psfrag{a2}{$a_2$}
\psfrag{a3}{$a_3$}
\psfrag{a4}{$a_4$}
\psfrag{aN}{$a_N$}
\psfrag{B1}{$b_1$}
\psfrag{B2}{$b_2$}
\psfrag{B3}{$b_3$}
\psfrag{B4}{$b_4$}
\psfrag{Bk}{$b_k$}
\psfrag{BN}{$b_N$}
\psfrag{hd}{$\cdots$}
\psfrag{vd}{$\vdots$}
\psfrag{_a31}{$a^3_1$}
\psfrag{_a32}{$a^3_2$}
\psfrag{_a21}{$a^2_1$}
\psfrag{_a22}{$a^2_2$}
\psfrag{_a23}{$a^2_3$}
\psfrag{_a24}{$a^2_4$}
\psfrag{_a11}{$a^1_1$}
\psfrag{_a12}{$a^1_2$}
\psfrag{_a13}{$a^1_3$}
\psfrag{_a01}{$a^0_1$}
\psfrag{_a02}{$a^0_2$}
\psfrag{C1}{$c_1$}
\psfrag{C2}{$c_2$}
\psfrag{C3}{$c_3$}
\psfrag{_1}{$1$}
\psfrag{AA0}{$a^0$}
\psfrag{AA1}{$a^1$}
\psfrag{AA2}{$a^2$}
\psfrag{AA3}{$a^3$}
\psfrag{BB0}{$b^0$}
\psfrag{BB1}{$b^1$}
\psfrag{BB2}{$b^2$}
\psfrag{BB3}{$b^3$}

\begin{document}

\title{Alternating signs of quiver coefficients}
\author{Anders Skovsted Buch}
\address{Matematisk Institut, Aarhus Universitet, Ny Munkegade, 8000
  {\AA}rhus C, Denmark}
\email{abuch@imf.au.dk}
\date{\today}
\maketitle

\section{Introduction}

Let $X$ be a non-singular algebraic variety and $E_0 \to E_1 \to
\cdots \to E_n$ a sequence of vector bundles and bundle maps over $X$.
A set of {\em rank conditions\/} for this sequence is a collection $r
= \{ r_{ij} \}$ of non-negative integers, for $0 \leq i < j \leq n$.
This data defines the quiver variety
\[ \Omega_r(E_\bull) = \{ x \in X \mid
   \rank(E_i(x) \to E_j(x)) \leq r_{ij} ~\forall i < j \} \,,
\]
which comes with a natural structure of subscheme of $X$, given as the
intersection of the zero sections of the maps $\bigwedge^{r_{ij}+1}
E_i \to \bigwedge^{r_{ij}+1} E_j$.  We demand that the rank conditions
can {\em occur\/} as the ranks over a point in $X$.  If we set $r_{ii}
= \rank(E_i)$, then this is equivalent to the conditions that $r_{ij}
\leq \min ( r_{i,j-1}, r_{i+1,j} )$ for all $0 \leq i < j \leq n$ and
$r_{ij} + r_{i-1,j+1} \geq r_{i-1,j} + r_{i,j+1}$ for all $0 < i \leq
j < n$.  In this case the expected codimension of the quiver variety
$\Omega_r(E_\bull)$ is the integer $d(r) = \sum_{i < j}
(r_{i,j-1}-r_{ij}) (r_{i+1,j}-r_{ij})$.

In joint work with Fulton \cite{buch.fulton:chern} we established a
formula for the cohomology (or Chow) class of the quiver variety when
this codimension is attained.  This was generalized in
\cite{buch:grothendieck} to the following formula for the structure
sheaf of a quiver variety in the Grothendieck ring $K(X)$ of algebraic
vector bundles on $X$:
\begin{equation} \label{E:orig_kqf}
   [\OO_{\Omega_r(E_\bull)}] \ = \ \sum_\mu
   c_\mu(r)\, G_{\mu_1}(E_1-E_0)\, G_{\mu_2}(E_2-E_1) \cdots
   G_{\mu_n}(E_n-E_{n-1}) \,.
\end{equation}
Here the sum is over finitely many sequences $\mu =
(\mu_1,\dots,\mu_n)$ of partitions $\mu_i$ such that the sum $\sum
|\mu_i|$ of the weights is greater than or equal to $d(r)$.  The
stable Grothendieck polynomials $G_{\mu_i}(E_i-E_{i-1})$ are defined
in section \ref{sec:grothpoly}.  The {\em quiver coefficients\/}
$c_\mu(r)$ appearing in this formula are integers which are uniquely
determined by the condition that (\ref{E:orig_kqf}) is true for all
varieties $X$ and bundle sequences $E_\bull$, together with the
condition that $c_\mu(r) = c_\mu(r+m)$ holds for all $m \in \N$, where
$r+m = \{ r_{ij} + m \}$ is the rank conditions obtained by adding the
integer $m$ to the original rank conditions.

A formula for quiver coefficients was also given in
\cite{buch.fulton:chern, buch:grothendieck}; in the case of
$K$-theory, this is based on the algebra of stable Grothendieck
polynomials constructed in \cite{buch:littlewood-richardson}.
Although the original formulas for quiver coefficients do not keep
track of their signs, it was conjectured that the cohomological quiver
coefficients (given by sequences $\mu$ such that $\sum |\mu_i| =
d(r)$) are non-negative, while the $K$-theoretic quiver coefficients
have signs that alternate with codimension, i.e.\ $(-1)^{\sum |\mu_i|
  - d(r)} c_\mu(r) \geq 0$.  Special cases of of these conjectures
have been proved in \cite{buch.fulton:chern, buch:on, buch:stanley,
  buch:grothendieck, buch.kresch.ea:schubert,
  buch.kresch.ea:grothendieck}.

In their recent paper \cite{knutson.miller.ea:four}, Knutson, Miller,
and Shimozono deliver a breakthrough within the theory of quiver
formulas, and prove at least two explicit combinatorial formulas for
the cohomological quiver coefficients, which show that these
coefficients are non-negative.  One of the important ideas in their
work is to reinterpret the {\em lace diagrams\/} of Abeasis and
Del-Fra \cite{abeasis.del-fra:degenerations} as sequences of partial
permutations.  This interpretation is explained by a Gr{\"o}bner
degeneration of a quiver variety in a matrix space into a union of
products of matrix Schubert varieties.  The {\em component formula\/}
of \cite{knutson.miller.ea:four} writes the cohomology class of a
quiver variety as a sum, over all `{\em minimal\/}' lace diagrams, of
the products of the Schubert polynomials for the corresponding partial
permutations.  The proof that quiver coefficients are non-negative is
obtained by proving a stable version of this component formula, where
the Schubert polynomials are replaced with Stanley symmetric
functions.  This is sufficient because Stanley symmetric functions are
known to be Schur positive \cite{edelman.greene:balanced,
  lascoux.schutzenberger:structure}.

Knutson, Miller, and Shimozono also prove a {\em ratio formula}, which
writes the class of a quiver variety as a quotient of two Schubert
polynomials.  This formula follows from a careful analysis of the
Zelevinsky map \cite{zelevinsky:two, lakshmibai.magyar:degeneracy},
and is in fact established for both cohomology and $K$-theory.  The
component formulas are proved using a combination of the Gr{\"o}bner
degeneration and the ratio formula, as well as a combinatorial study
of a double version of the ratio formula.  In particular, it is proved
that a limit of the double ratio formula agrees with the {\em double
  quiver functions\/} introduced in \cite{buch.fulton:chern,
  buch:grothendieck} and named in \cite{knutson.miller.ea:four}.  The
authors of \cite{knutson.miller.ea:four} have informed us that they
can generalize their methods to work in $K$-theory, although,
according to their own description, this approach is rather
complicated.\footnote{This work is now available in the paper
  \cite{miller:alternating}, which proves a $K$-theoretic component
  formula and a double version of it, but only in the stable case and
  for large rank conditions.}

In this paper we give simpler proofs of the above mentioned formulas,
using methods that work equally well in $K$-theory.  In particular, we
prove that the $K$-theoretic quiver coefficients have alternating
signs, and we derive an explicit combinatorial formula for these
coefficients.  Starting from the ratio formula, we give combinatorial
proofs of $K$-theoretic generalizations of the component formulas,
where the Schubert polynomials and Stanley symmetric functions are
replaced with ordinary and stable Grothendieck polynomials.  These
formulas are given in terms of sequences of partial permutations,
which we call {\em KMS-factorizations\/} of the {\em Zelevinsky
  permutation\/} defined in \cite{knutson.miller.ea:four}.  To
conclude that $K$-theoretic quiver coefficients have alternating
signs, we use Lascoux's result that stable Grothendieck polynomials
are linear combinations with alternating signs of Stable Grothendieck
polynomials for partitions \cite{lascoux:transition}.  To make our
paper self-contained, we also give a short proof of the ratio formula.

The {\em factor sequences conjecture\/} of \cite{buch.fulton:chern}
states that cohomological quiver coefficients count the number of
sequences of semistandard Young tableaux which can be generated by a
sequence of factorizations and multiplications of chosen tableaux
arranged in a {\em tableaux diagram}.  A special case of this
conjecture, corresponding to a particular choice of tableaux diagram,
was proved in \cite{knutson.miller.ea:four}.  However, so far there
has been no progress in generalizing this conjecture to $K$-theory.
In this paper we close this gap by showing that KMS-factorizations can
be defined by the same algorithm as defines factor sequences, except
that the tableau diagram is replaced with a diagram of permutations,
and the plactic product of tableaux is replaced with multiplication of
permutations in the degenerate Hecke algebra.  A tableau-based version
of this type of factor sequences also exists; we briefly outline this
for cohomological quiver coefficients, and refer to \cite{bkty:stable}
for details and the general case.

We remark that it was already known that the cohomological component
formulas can be derived combinatorially from the ratio formula, by
using a simplification of Yong \cite{yong:embedding}.  However, Yong's
method still requires the analysis of the double ratio formula and its
limits from \cite{knutson.miller.ea:four}.  The approach presented
here simplifies things further by working only with the single ratio
formula, by applying Fomin and Kirillov's construction of Grothendieck
polynomials based on solutions to Yang-Baxter equations
\cite{fomin.kirillov:yang-baxter, fomin.kirillov:grothendieck}, and by
observing that the stable component formula follows easily from the
non-stable component formula.

Other simple proofs of the cohomological component formulas have also
surfaced.  For example, they can be deduced very easily from the Thom
polynomial theory developed by Feh\'er and Rim\'anyi
\cite{feher.rimanyi:classes}, or deduced directly from the above
mentioned Gr{\"o}bner degeneration with a symmetry argument.  This is
explained in \cite{buch.feher.ea:positivity}.  While attempts to
generalize these methods to $K$-theory have not been successful, they
might hold more promise for quiver varieties of other types (see
\cite{feher.rimanyi:classes}).

Some of the results proved in \cite{knutson.miller.ea:four} imply that
the cohomological double ratio formula for large rank conditions
satisfies nice properties, including {\em multi supersymmetry\/} and a
double version of the component formula.  In the last section of this
paper, we establish these properties for the $K$-theoretic double
ratio formula given by arbitrary rank conditions.  In particular, we
prove a conjecture from \cite{knutson.miller.ea:four} stating that the
double ratio formula satisfies a {\em rank stability\/} property.
This conjecture is equivalent to the statement that the polynomials
defined by the double ratio formula are specializations of the
original quiver formulas \cite{buch.fulton:chern, buch:grothendieck}.
Even though the double ratio formula is not needed for the proof of
alternating signs of quiver coefficients given in this paper, its
multi supersymmetry property has some nice applications.  For example,
this property was used in \cite{knutson.miller.ea:four} to prove the
above mentioned case of the factor sequences conjecture.  The multi
supersymmetry property also implies that general quiver coefficients
are special cases of the coefficients studied in
\cite{buch.kresch.ea:grothendieck}.  In fact, quiver coefficients can
be realized as Schubert structure constants on flag varieties
\cite{bergeron.sottile:schubert, lenart.robinson.ea:grothendieck,
  buch.sottile.ea:quiver}.

This paper is organized as follows.  In section \ref{sec:grothpoly} we
explain Fomin and Kirillov's construction of Grothendieck polynomials,
which is a key ingredient in the combinatorial parts of this paper.
Section \ref{sec:zel} gives a new construction of the Zelevinsky
permutation which is required for our proof that KMS-factorizations
can be viewed as factor sequences.  Section \ref{sec:ratio} contains
the proof of the ratio formula.  In section~\ref{sec:restrict} we
prove a formula for double Grothendieck polynomials applied to certain
rearrangements of the same set of variables, which in section
\ref{sec:component} is used to derive the non-stable component formula
from the ratio formula.  In section \ref{sec:facseq} we establish the
factor sequences definition of KMS-factorizations and discuss its
consequences.  As a corollary we obtain a rank stability property for
KMS-factorizations, which in section \ref{sec:altsigns} is used to
derive the stable component formula and deduce that quiver
coefficients have alternating signs.  Section \ref{sec:double} finally
proves the above mentioned properties of the double ratio formula.

We are very grateful to Rich\'ard Rim\'anyi for discussions at the
Banach Institute in Warsaw, which led to our observation that the
stable component formula can be deduced from the non-stable formula,
and which triggered our search for other simplifications to
\cite{knutson.miller.ea:four}.  We also thank Feh{\'e}r, Kresch,
Sottile, Tamvakis, and Yong for inspiring collaboration on related
papers, and Fulton, Miller, and Sottile for helpful comments to our
paper.  Finally, we thank Martin Guest, Anatol Kirillov, and the
Research Institute for Mathematical Sciences in Kyoto for their
hospitality while this paper was written.


\section{Grothendieck polynomials}
\label{sec:grothpoly}

The {\em degenerate Hecke algebra\/} $\Hecke$ over a commutative ring
$R$ is the free $R$-algebra generated by symbols $s_1, s_2, \dots$,
modulo the relations
\begin{eqnarray*}
  s_i s_j &=& s_j s_i  \text{\ \ \ \ \ \ if $|i-j| \geq 2$} \\
  s_i s_{i+1} s_i &=& s_{i+1} s_i s_{i+1} \\
  s_i^2 &=& - s_i \,.
\end{eqnarray*}
In this paper, $R$ will be a ring of Laurent polynomials.  The algebra
$\Hecke$ is a free $R$-module with a basis of permutations
corresponding to reduced expressions in the generators.

Given permutations $u_1, u_2, \dots, u_n$, the product $u_1 \cdot u_2
\cdots u_n$ in $\Hecke$ of these permutations is equal to plus or
minus a single permutation $w$.  We will call this permutation $w$ for
the {\em absolute Hecke product\/} of the $u_i$.  Notice that the
descent positions of $w$, i.e.\ the indices $i$ for which $w(i) >
w(i+1)$, include the descent positions of $u_n$, while the descent
positions of $u_1^{-1}$ are also descent positions of $w^{-1}$.
Notice also that if $u_1 u_2 \cdots u_n$ is a {\em reduced product\/}
of permutations in the sense that $\ell(u_1 u_2 \cdots u_n) = \sum
\ell(u_i)$, then the absolute Hecke product $w$ agrees with the usual
product of permutations.

We will need Fomin and Kirillov's construction
\cite{fomin.kirillov:yang-baxter, fomin.kirillov:grothendieck} of the
Grothendieck (Laurent) polynomials $\Groth_w(a;b)$ of Lascoux and
Sch{\"u}tzenberger \cite{lascoux.schutzenberger:structure,
  lascoux:anneau}.  Consider a diagram $\D$ of strings going from a
top horizontal border to a left vertical border.  Each string must be
composed of line segments, each of which is labeled with a variable
and has a direction between due south and due west.  Furthermore,
strings may only cross each other transversally, at inner points of
the line segments.  Of particular importance is the following diagram
$\D_N$, which contains only horizontal and vertical line segments.
\[ \D_N \ \ = \ \ \ \ \raisebox{-40pt}{\pic{60}{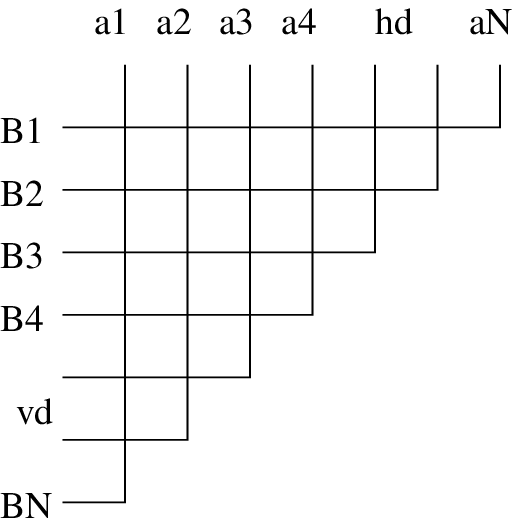}} \]

Let $C(\D)$ denote the set of crossing positions in the diagram $\D$.
For each $P \in C(\D)$ we set $h(P) = 1 - \frac{y}{x}$ where $x$ is
the label of the line through $P$ with the highest slope (within the
range $[0,+\infty]$), and $y$ is the label of the line with the lowest
slope.  We also let $\nu(P)$ be one plus the number of strings in $\D$
passing north-west of $P$.  We then define the {\em FK-product\/}
$\Groth(\D)$ as the product in $\Hecke$ of the factors $(1 + h(P)\,
s_{\nu(P)})$ for all $P \in C(\D)$.  These factors should be
multiplied from south-west to north-east, in any order so that each
crossing position $P$ comes before all other crossing positions in the
quadrangle between the lines going due north and due east from $P$.

For the diagram $\D_N$, we identify the crossing point of the
horizontal line labeled $b_p$ and the vertical line labeled $a_q$,
with the point $(p,q) \in \N \times \N$.  We then have $\nu(p,q) =
p+q-1$, and the FK-product of $\D_N$ is given by
\[ \Groth(\D_N) = \prod_{q=1}^{N-1} \prod_{p = N-q}^1 
   \left( 1 + \left(1-\frac{b_p}{a_q}\right) s_{p+q-1} \right) 
   \in \Hecke \,.
\]
We need the following theorem which is proved in
\cite{fomin.kirillov:grothendieck} (modulo the change of variables
$x_i = 1 - a_i^{-1}$ and $y_i = 1 - b_i$; see Thm.~2.3 and the remark
on page 7 of loc.\ cit.)

\begin{thm}[Fomin and Kirillov] \label{T:stddiag}
  In $\Hecke$ we have the identity
\[
  \Groth(\D_N) = \sum_{w\in S_N} \Groth_w(a;b) \cdot w
\]
where $\Groth_w(a;b)$ is the double Grothendieck polynomial for $w$.
\end{thm}

Suppose $D$ is a subset of the crossing positions $C(\D)$ of a diagram
$\D$.  We let $w(D)$ be the absolute Hecke product of the simple
reflections $s_{\nu(P)}$ for $P \in D$, in south-west to north-east
order as above.  We say that $D$ is an {\em FK-graph\/} for this
permutation $w(D)$, and that $D$ is {\em reduced\/} if $|D| =
\ell(w(D))$.  We can picture an FK-graph $D$ by replacing the crossing
positions of $\D$ which belong to this graph with the symbol
``\,$\cross$\,'', while the remaining crossing positions are replaced
with the symbol ``\,$\avoid$\,''.  Notice that if $D$ is reduced, then
the string entering the resulting diagram at column $i$ at the top
will exit at row $w(D)(i)$ at the left hand side.
Notice also that any FK-graph $D$ contains a reduced FK-graph
$D' \subset D$ such that $w(D') = w(D)$.  In fact, $D'$ can be found
by simply skipping the points $P \in D$ for which $s_{\nu(P)}$ does
not increase the length when the product $w(D)$ is formed.


When no diagram $\D$ is explicitly mentioned, an FK-graph will always
be relative to a diagram $\D_N$, so it is a finite subset of $\N
\times \N$.  Such FK-graphs are called {\em pipe dreams\/} in
\cite{knutson.miller:grobner}, and a reduced pipe dream is the same as
an {\em RC-graph\/} \cite{bergeron.billey:rc-graphs}.  For example,
the pipe dream
\[ D \ \ = \ \ \raisebox{-18pt}{\pic{14}{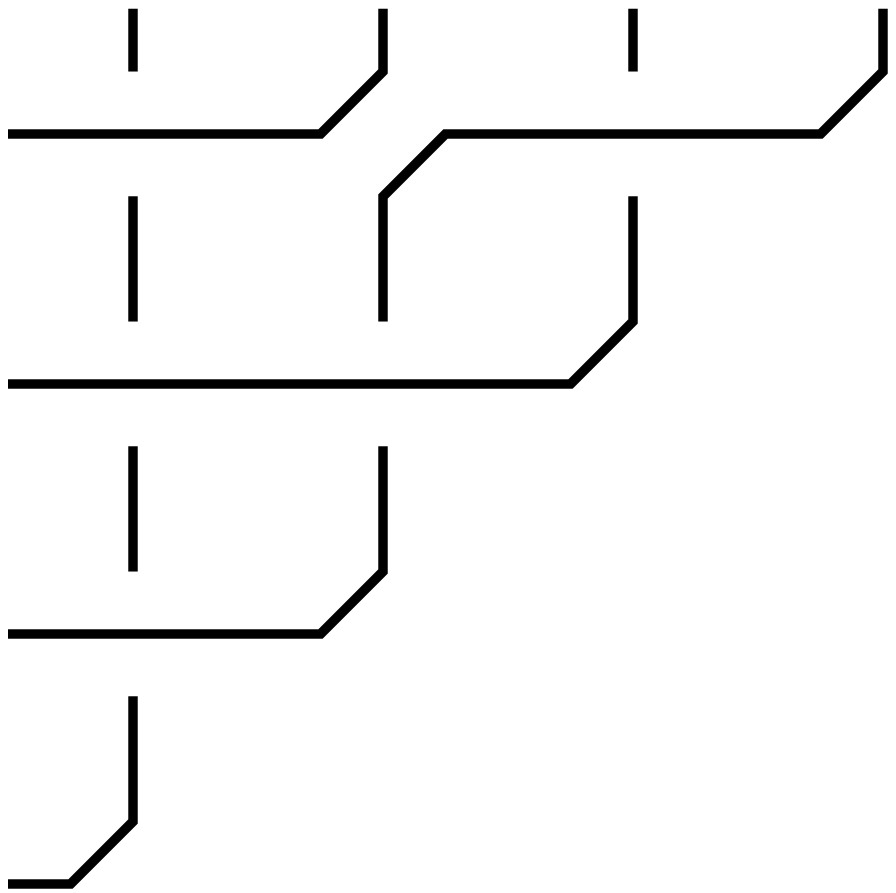}} \]
is an FK-graph for the permutation $w(D) = s_3 s_2 s_1 s_3 = 4132$.

It follows from the definitions that the coefficient of a permutation
$w$ in the FK-product of a diagram $\D$ is equal to
\begin{equation} \label{E:fkprd}
  \sum_{D \subset C(\D) \,,\, w(D) = w} (-1)^{|D|-\ell(w)}\,
  \prod_{P \in D} h(P) \,.
\end{equation}
\refthm{T:stddiag} therefore has the following corollary.  See also
\cite{knutson.miller:subword} for an alternative proof and
\cite{bergeron.billey:rc-graphs} for the case of Schubert polynomials.

\begin{cor} \label{C:fkgraph}
For any permutation $w$ we have
\[ \Groth_w(a;b) = \sum_{w(D)=w} (-1)^{|D|-\ell(w)} 
   \prod_{(p,q) \in D} \left( 1 - \frac{b_p}{a_q} \right) 
\]
where the sum is over all FK-graphs $D \subset \N \times \N$ for
$w$.
\end{cor}

Let us remark that \refthm{T:stddiag} is more flexible than its
corollary, as amply demonstrated in \cite{fomin.kirillov:yang-baxter}.
The point is that many operations can be performed on a diagram $\D$
without changing the corresponding FK-product.  We will write $\D
\approx \D'$ if $\Groth(\D) = \Groth(\D')$.  The two key examples of
this are:
\begin{equation} \label{E:key}
  \raisebox{-18pt}{\pic{50}{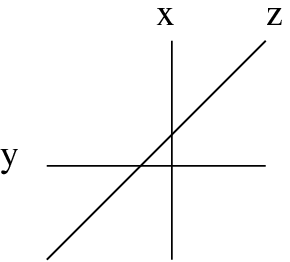}} \ \approx \ \ 
  \raisebox{-18pt}{\pic{50}{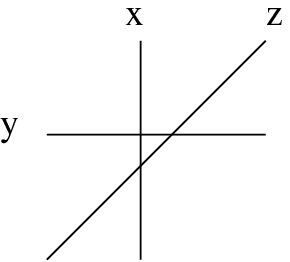}}  
  \ \ \ \ \ \ \text{and} \ \ \ \ \ \ 
  \raisebox{-18pt}{\pic{50}{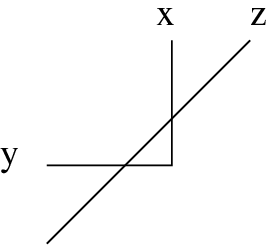}} \ \approx \ \  
  \raisebox{-18pt}{\pic{50}{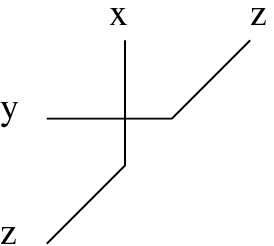}} 
\end{equation}
Notice that when $x=y$, the last diagram is also equivalent to a north
to west hook labeled $x$ together with a disjoint north-east to
south-west line labeled $z$.  Although most diagrams in this paper
contain only horizontal and vertical line segments, the availability
of slanted lines often makes it more natural to manipulate these
diagrams using the rules of (\ref{E:key}).

Lascoux and Sch\"utzenberger's original definition of Grothendieck
polynomials says that $\Groth_{w_0}(a;b) = \prod_{p+q \leq N}
(1-b_p/a_q)$ where $w_0 \in S_N$ is the longest permutation, and that
$(a_i-a_{i+1})\Groth_w(a;b) = a_i \Groth_{w s_i}(a;b) - a_{i+1}
\Groth_{w s_i}(a_{s_i};b)$ when $w(i)<w(i+1)$.  Here we let $a_v$
denote the sequence of variables $a_{v(1)},\dots,a_{v(N)}$ for any
permutation $v \in S_N$.  It follows that if $i$ is not a descent
position for $w$, then $\Groth_w(a;b)$ is symmetric in the variables
$a_i$ and $a_{i+1}$.  In particular, if $k$ is the last descent
position for $w$, then the variables $a_i$ for $i > k$ do not occur in
$\Groth_w(a;b)$.  A similar relationship holds between the descent
positions of $w^{-1}$ and the variables $b_i$.  We need the identity
\begin{equation} \label{E:conjinv}
  \Groth_w(a_{v(1)},\dots,a_{v(N)};a_1,\dots,a_N) = 
  \Groth_{w_0 w^{-1} w_0}(a_N,\dots,a_1; a_{v(N)},\dots,a_{v(1)}) \,.
\end{equation}
In fact, using that $\Groth_w(a;b) =
\Groth_{w^{-1}}(b_1^{-1},\dots,b_N^{-1};a_1^{-1},\dots,a_N^{-1})$ and
the above definition, equation (\ref{E:conjinv}) follows by descending
induction on $\ell(w)$ from the calculation
\[\begin{split}
  & (a_{v(i)}-a_{v(i+1)}) \Groth_w(a_v;a) = a_{v(i)} 
    \Groth_{w s_i}(a_v;a) - a_{v(i+1)} \Groth_{w s_i}(a_{v s_i};a) \\
  &\hspace{50pt} = 
    a_{v(i)} \Groth_{w_0 s_i w^{-1} w_0}(a_{w_0}; a_{v w_0}) -
    a_{v(i+1)}\Groth_{w_0 s_i w^{-1} w_0}(a_{w_0}; a_{v s_i w_0}) \\
  &\hspace{50pt} = 
    (a_{v(i)}-a_{v(i+1)}) \Groth_{w_0 w^{-1} w_0}(a_{w_0};a_{v w_0}) \,.
\end{split}\]
The identity (\ref{E:conjinv}) also follows from the results in
\cite{lascoux.schutzenberger:decompositions}.

\ignore{
It follows from Lascoux and Sch\"utzenberger's original definition of
Grothendieck polynomials that, if $i$ is not a descent position for
$w$, then $\Groth_w(a;b)$ is symmetric in the variables $a_i$ and
$a_{i+1}$.  In particular, if $k$ is the last descent position for
$w$, then the variables $a_i$ for $i > k$ do not occur in
$\Groth_w(a;b)$.  A similar relationship holds between the descent
positions of $w^{-1}$ and the variables $b_i$.
}

We also need the {\em stable double Grothendieck polynomials\/}
$G_w(a; b)$ of Fomin and Kirillov \cite{fomin.kirillov:grothendieck}.
These polynomials are characterized by the property that
\[ G_w(a_1,\dots,a_q; b_1,\dots,b_p) = \Groth_{1^m \times
  w}(a_1,\dots,a_q,1,\dots,1 \,;\, b_1,\dots,b_p,1,\dots,1)
\]
for all $m \geq \max(p,q)$.  Here the permutation $1^m \times w$ is
the identity on $\{1,\dots,m\}$ while it maps $j$ to $w(j-m)+m$ for $j
> m$.  If $\lambda = (\lambda_1 \geq \dots \geq \lambda_k)$ is a
partition, we set $G_\lambda(a;b) = G_{w_\lambda}(a;b)$ where
$w_\lambda$ is the Grassmannian permutation for $\lambda$, defined by
$w_\lambda(i) = i + \lambda_{k+1-i}$ for $1 \leq i \leq k$ and
$w_\lambda(i) < w_\lambda(i+1)$ for $i \neq k$.  

It is proved in \cite{buch:littlewood-richardson} that any stable
Grothendieck polynomial $G_w(a;b)$ can be written as an integral
linear combination
\begin{equation} \label{E:stabcoef}
  G_w(a;b) = \sum_\lambda c_{w,\lambda} \, G_\lambda(a;b)
\end{equation}
of stable Grothendieck polynomials for partitions.  Lascoux has proved
\cite[Thm.~4]{lascoux:transition} an explicit combinatorial formula
for the coefficients $c_{w,\lambda}$ in this expansion, which shows
that they have alternating signs, {i.e.\ }$(-1)^{|\lambda|-\ell(w)}\,
c_{w,\lambda} \geq 0$.  (See also the reformulation of Lascoux's
formula in \cite[Thm.~3]{buch.kresch.ea:grothendieck}.)

Given vector bundles $F = L_1 \oplus \dots \oplus L_p$ and $H = M_1
\oplus \dots \oplus M_q$ over $X$ which are direct sums of line
bundles, we write $G_\lambda(H-F) = G_\lambda(M_1,\dots,M_q ;
L_1,\dots,L_p) \in K(X)$.  By the symmetry of $G_\lambda$, this is a
polynomial in the exterior powers of $F$ and $H^\vee$, so the notation
$G_\lambda(H-F)$ also makes sense for bundles which are not sums of
line bundles.  This explains the notation used in the quiver formula
(\ref{E:orig_kqf}).

\ignore{
We shall need the following property of Grothendieck polynomials.  If
$v, w \in S_N$ are permutations and $w_0 \in S_N$ is the longest
permutation, then we have
\begin{equation} \label{E:conjinv}
  \Groth_w(a_{v(1)},\dots,a_{v(N)};a_1,\dots,a_N) = 
  \Groth_{w_0 w^{-1} w_0}(a_N,\dots,a_1; a_{v(N)},\dots,a_{v(1)}) \,.
\end{equation}
If we let $a_v$ denote the sequence of variables
$a_{v(1)},\dots,a_{v(N)}$, then Lascoux and Sch\"utzenberger's
original definition of Grothendieck polynomials states that
$(a_i-a_{i+1})\Groth_w(a;b) = a_i \Groth_{w s_i}(a;b) - a_{i+1}
\Groth_{w s_i}(a_{s_i};b)$ when $w(i)<w(i+1)$.  Using this and the
identity $\Groth_w(a;b) =
\Groth_{w^{-1}}(b_1^{-1},\dots,b_N^{-1};a_1^{-1},\dots,a_N^{-1})$ one
checks that
\[ (a_{v(i)}-a_{v(i+1)}) \Groth_w(a_v;a) = a_{v(i)} \Groth_{w
   s_i}(a_v;a) - a_{v(i+1)} \Groth_{w s_i}(a_{v s_i};a)
\]
and
\begin{multline*}
  (a_{v(i)}-a_{v(i+1)}) \Groth_{w_0 w^{-1} w_0}(a_{w_0};a_{v w_0}) = \\
  a_{v(i)} \Groth_{w_0 s_i w^{-1} w_0}(a_{w_0}; a_{v w_0}) -
  a_{v(i+1)}\Groth_{w_0 s_i w^{-1} w_0}(a_{w_0}; a_{v s_i w_0}) \,.
\end{multline*}
The identity (\ref{E:conjinv}) follows from these identities by
descending induction on $\ell(w)$.
}


\section{The Zelevinsky permutation}
\label{sec:zel}

In this section we give a new construction of the Zelevinsky
permutation of \cite{knutson.miller.ea:four}, which is needed for our
proof of the $K$-theoretic analogue of the factor sequences
conjecture.  To be precise, we construct the conjugate of the
Zelevinsky permutation, which turns out to have a simpler relationship
with the geometry of quiver varieties and KMS-factorizations.
However, the Zelevinsky permutation itself is necessary to obtain nice
combinatorial properties of the ratio formula.

Extend the set of rank conditions $r$ by setting $r_{ij} = e_j +
e_{j+1} + \cdots + e_i$ for $j \leq i$, where $(e_0,\dots,e_n) =
(r_{00},\dots,r_{nn})$ is the dimension vector corresponding to $r$.
For $i < 0$ or $j > n$ we set $r_{ij}=0$, and we set $N = r_{n0} =
e_0+\dots+e_n$.  For each $0 \leq i < n$ and $0 < j \leq n$ we define
a permutation $W_{ij} \in S_{r_{i+1,j-1}}$ by the expression
\[ W_{ij}(p) = \begin{cases}
   p + r_{i,j-1} - r_{ij} & \text{if $r_{ij} < p \leq r_{i+1,j}$} \\
   p - r_{i+1,j} + r_{ij} & \text{if $r_{i+1,j} < p \leq 
                                  r_{i+1,j} + r_{i,j-1} - r_{ij}$} \\
   p & \text{otherwise.}
\end{cases}\]
When $i < j$, this is the Grassmannian permutation for the rectangular
partition $R_{ij}$ from \cite{buch.fulton:chern}, with descent at
position $r_{i+1,j}$.  Now define the {\em conjugate Zelevinsky
  permutation\/} $z(r) \in S_N$ for the rank conditions $r$ to be the
south-west to north-east product of the matrix of permutations
\begin{equation} \label{E:permmat}
\begin{matrix}
  W_{n-1,1} & W_{n-1,2} & \cdots & W_{n-1,n} \\
  W_{n-2,1} & W_{n-2,2} & \cdots & W_{n-2,n} \\
  \vdots    & \vdots    & \ddots & \vdots \\
  W_{0,1}   & W_{0,2}   & \cdots & W_{0,n} 
\end{matrix} \ ,
\end{equation}
that is $z(r) = \prod_{j=1}^n \prod_{i=0}^{n-1} W_{ij}$.  Notice that
a `south-west to north-east product' makes sense because $W_{ij}$
commutes with $W_{i'j'}$ whenever $i < i'$ and $j' < j$.  Notice also
that $W_{ij}$ is the conjugate Zelevinsky permutation for the set of
rank conditions consisting of the integers $r_{i,j-1}$, $r_{i+1,j}$,
and $r_{ij}$.  Let $w_0^{(N)}$ denote the longest permutation in
$S_N$.  We will call the permutation $\wh z(r) = w_0^{(N)}
z(r)^{-1} w_0^{(N)}$ for the {\em Zelevinsky permutation\/} for $r$,
although the original definition in \cite{knutson.miller.ea:four}
assigns this name to the inverse of $\wh z(r)$.  The action of $z(r)$
can be described explicitly as follows.

\begin{lemma} \label{L:zeldesc}
  Given any integer $1 \leq p \leq N$, there are unique integers $i$
  and $j$ with $0 \leq i \leq j+1 \leq n+1$, such that
  $r_{n,j+1}+r_{i-1,j}-r_{i-1,j+1} < p \leq
  r_{n,j+1}+r_{i,j}-r_{i,j+1}$.  We then have $z(r)(p) = p - r_{n,j+1}
  + r_{i,j+1} + r_{i-1,0} - r_{i-1,j}$.
\end{lemma}
\begin{proof}
  To find $i$ and $j$, one first chooses $j$ such that $r_{n,j+1} < p
  \leq r_{nj}$, after which $i$ is uniquely determined.  Now, when
  $z(r)$ is applied to $p$, the permutation $W_{k,j+1}$ subtracts
  $r_{k+1,j+1} - r_{k,j+1}$ from its argument for $k=n-1,n-2,\dots,i$,
  after which $W_{i-1,k}$ adds $r_{i-1,k-1} - r_{i-1,k}$ to its
  argument for $k = j,j-1,\dots,1$.  All other factors of $z(r)$
  preserve their argument.
\end{proof}

The relation of the Zelevinsky permutation with quiver varieties is
based on the following lemma, which is equivalent to
\cite[Prop.~1.6]{knutson.miller.ea:four}.

\begin{lemma} \label{L:zelprop}
  The conjugate Zelevinsky permutation $z(r)$ is the unique
  permutation in $S_N$ such that
  
  (i) all descent positions of $z(r)$ are contained in the set $\{
  r_{nj} \mid 0 < j \leq n \}$,
  
  (ii) all descent positions of $z(r)^{-1}$ are contained in $\{ r_{i0}
  \mid 0 \leq i < n \}$, and
  
  (iii) for all $0 \leq i,j \leq n$, there are $r_{ij}$ integers $p
  \leq r_{nj}$ for which $z(r)(p) \leq r_{i0}$.
\end{lemma}
\begin{proof}
  Part (i) follows directly from \reflemma{L:zeldesc}, and (ii)
  follows by observing that $z(r)^{-1}$ is the conjugate Zelevinsky
  permutation for the mirrored rank conditions $r'$ given by $r'_{ij}
  = r_{n-j,n-i}$.  Notice that if $p,i,j$ are chosen as in
  \reflemma{L:zeldesc} then $r_{i-1,0} < z(r)(p) \leq r_{i,0}$.
  Therefore there are exactly $r_{i,j} - r_{i,j+1} - r_{i-1,j} +
  r_{i-1,j+1}$ integers $p$ such that $r_{n,j+1} < p \leq r_{n,j}$ and
  $r_{i-1,0} < z(r)(p) \leq r_{i,0}$.  This proves (iii) since $r_{ij}
  = \sum_{k=0}^i \sum_{l=j}^n (r_{kl} - r_{k,l+1} - r_{k-1,l} +
  r_{k-1,l+1})$.  The uniqueness statement is not needed in this
  paper, and its easy proof is left to the reader.
\end{proof}

\reflemma{L:zelprop} implies that the length of $z(r)$ is given by
\begin{equation}
  \ell(z(r)) = \sum_{0\leq i<n, 0<j\leq n}
  (r_{i,j-1}-r_{ij})(r_{i+1,j}-r_{ij}) \,.
\end{equation}
In fact, there are $(r_{i,j-1}-r_{ij})(r_{i+1,j}-r_{ij})$ pairs
$(p,q)$ for which $p \leq r_{nj} < q \leq r_{n,j-1}$ and $z(r)(q) \leq
r_{i0} < z(r)(p) \leq r_{i+1,0}$.  In particular, $z(r)$ is a reduced
product of the permutations $W_{ij}$ of the matrix (\ref{E:permmat}).

We let $z(e) = z(r^e)$ be the conjugate Zelevinsky permutation for the
maximal rank conditions $r^e$ given by $r^e_{ij} = \min \{e_i,
e_{i+1},\dots,e_j\}$ for $i \leq j$.  Similarly we write $\wh z(e) =
w_0^{(N)} z(e)^{-1} w_0^{(N)}$.  The inverse of this permutation is
called $v(Hom)$ in \cite{knutson.miller.ea:four}.

It follows from \reflemma{L:zeldesc} that for all $i$ and $p > r_{ni}$
we have $z(r)(p) \leq r_{i0}$.  This implies that for $q \leq
r_{n,i+1}$ we have $\wh z(r)(q) = N+1-z(r')(N+1-q) > r_{i-1,0}$, so
any FK-graph for $\wh z(r)$ must contain the set $D_e =
\bigcup_{i=1}^{n-1} [1,r_{i-1,0}] \times [1,r_{n,i+1}]$.  Since $|D_e|
= \sum_{j-i \geq 2} e_i e_j = \ell(\wh z(e))$, it follows that $D_e$
is the unique FK-graph for $\wh z(e)$.  \refcor{C:fkgraph} therefore
implies that the polynomial $\Groth_{\wh z(e)}(a;b) = \prod_{(p,q)\in
  D_e} (1 - \frac{b_p}{a_q})$ divides $\Groth_{\wh z(r)}(a;b)$ (cf.\ 
\cite[\S 5.1]{knutson.miller.ea:four}.)


\section{The ratio formula}
\label{sec:ratio}

In this section we give a coordinate free proof of the ratio formula
\cite[Thm.~2.7]{knutson.miller.ea:four}.  The underlying geometry is
similar to the original proof, but becomes slightly simpler by working
with the conjugate Zelevinsky permutation.

Suppose $F_1 \subset F_2 \subset \dots \subset F_N \to H_N \tto \cdots
\tto H_2 \tto H_1$ is a flag of vector bundles over $X$ with a
morphism to a dual flag, such that $\rank(F_i) = \rank(H_i) = i$.  For
a permutation $w \in S_N$, Fulton \cite{fulton:flags} defines the
degeneracy locus
\[ \Omega_w = \{ x \in X \mid \rank(F_q(x) \to H_p(x)) \leq r_w(p,q)
   ~\forall p,q \} 
\]
where $r_w(p,q) = \# \{ k \leq p \mid w(k) \leq q \}$.  This locus
does not depend on the bundles $H_p$ for which $w(p) < w(p+1)$, or on
the bundles $F_q$ such that $w^{-1}(q) < w^{-1}(q+1)$.  The expected
codimension of the locus is the length $\ell(w)$.  We need the
following formula for the Grothendieck class of $\Omega_w$, which
generalizes Fulton's formula for its cohomology class
\cite{fulton:flags}.  The $K$-theory formula was proved in
\cite[Thm.~2.1]{buch:grothendieck} as an application of
\cite[Thm.~3]{fulton.lascoux:pieri}.  An equivalent statement was
given in \cite[Thm.~A]{knutson.miller:grobner}.

\begin{thm} \label{T:fulton}
If the codimension of $\Omega_w$ in $X$ equals $\ell(w)$, then 
\[ [\OO_{\Omega_w}] = \Groth_w(L_1,\dots,L_N; M_1,\dots,M_N) \ \in K(X) \]
where $L_i = \ker(H_i \to H_{i-1})$ and $M_i = F_i/F_{i-1}$.
\end{thm}

It follows from \reflemma{L:zelprop} that the Grothendieck polynomial
$\Groth_{\wh z(r)}(a;b)$ is symmetric in each interval of variables
$\{ a_p : r_{n,i+1} < p \leq r_{n,i}\}$ and $\{b_p : r_{i-1,0} < p
\leq r_{i,0}\}$.  We therefore allow the vector bundle $E_i$ to be
substituted for these intervals, and interpret the result as if
imaginative line bundle summands of $E_i$ had been inserted.  We use
the same notation in the quotient of polynomials $\Groth_{\wh
  z(r)}(a;b)/\Groth_{\wh z(e)}(a;b)$.

\begin{thm}[Knutson, Miller, Shimozono] \label{T:ratio}
  If the codimension of $\Omega_r(E_\bull)$ in $X$ is equal to $d(r)$,
  then the class of its structure sheaf in $K(X)$ is given by
\[ [\OO_{\Omega_r(E_\bull)}] \ = \  
  \frac{\Groth_{\wh z(r)}(E_n,\dots,E_0 \,;\, E_0,\dots,E_n)}
  {\Groth_{\wh z(e)}(E_n,\dots,E_0 \,;\,E_0,\dots,E_n)} \,.
\]
\end{thm}
\begin{proof}
  By replacing $X$ with $\bigoplus_{j-i\geq 2} \Hom(E_{i-1},E_i)$, we
  may assume that there are general maps $\phi_{ij} : E_i \to E_j$ for
  all $i < j$, such that each map $\phi_{i-1,i}$ comes from the given
  bundle sequence $E_\bull$.  We may then construct the bundle
  sequence
%
\[ E_0 \subset E_0\oplus E_1 \subset \dots \subset
   E_0\oplus\dots\oplus E_n \xrightarrow{\phi}
   E_0\oplus\dots\oplus E_n \tto \dots \tto 
   E_{n-1}\oplus E_n \tto E_n
\]
where the middle map $\phi$ is composed of maps $E_i \to E_j$ which is
$\phi_{ij}$ for $i < j$, the identity for $i=j$, and zero for $i > j$.
All other maps are the obvious embeddings or projections.  It follows
from \reflemma{L:zelprop} that Fulton's locus $\Omega_{z(r)}$ for this
sequence consists of the points in $X$ where each composed map
$E_0\oplus\dots\oplus E_i \to E_j\oplus\dots\oplus E_n$ has rank at
most $r_{ij}$.

We claim that $\Omega_{z(r)}$ is the intersection of the quiver
variety $\Omega_r(E_\bull)$ with the subset of $X$ where each map
$\phi_{ij}$ is equal to the composition $\phi_{j-1,j}
\phi_{j-2,j-1}\cdots \phi_{i,i+1}$.  To see this, notice that the
condition that the rank of the map $E_0\oplus\dots\oplus E_k \to
E_k\oplus\dots\oplus E_n$ is at most $e_k$ is equivalent to demanding
that for all $i < k < j$ we have $\phi_{ij} = \phi_{kj} \phi_{ik}$.
Given that this holds for all $k$, we obtain that for each $i<j$ the
map $E_0\oplus\dots\oplus E_i \to E_j\oplus\dots\oplus E_n$ can be
factored as a surjection followed by $\phi_{ij}$ followed by an
injection, so the rank condition on this map is equivalent to
$\rank(\phi_{ij}) \leq r_{ij}$.  In particular, $\Omega_{z(e)}$ is the
locus where $\phi_{ij} = \phi_{j-1,j} \cdots \phi_{i,i+1}$ for all $i
< j$, and we have $\Omega_{z(r)} = \Omega_r(E_\bull) \cap
\Omega_{z(e)}$ as a scheme theoretic intersection.

The above description also shows that the codimension of
$\Omega_{z(r)}$ in $X$ is $\ell(z(r)) = d(r) + \ell(z(e))$, so by
\refthm{T:fulton} and (\ref{E:conjinv}) we have $[\OO_{\Omega_{z(r)}}]
= \Groth_{z(r)}(\breve E_\bull; E_\bull) = \Groth_{\wh z(r)}(\breve
E_\bull; E_\bull)$, where we write $\breve E_\bull$ for the reversed
sequence $E_n,E_{n-1},\dots,E_0$.  Similarly we have
$[\OO_{\Omega_{z(e)}}] = \Groth_{\wh z(e)}(\breve E_\bull; E_\bull)$.
Since all of the degeneracy loci are Cohen-Macaulay
\cite{lakshmibai.magyar:degeneracy} we obtain the identity
\begin{equation} \label{E:ratio}
  \Groth_{\wh z(r)}(\breve E_\bull; E_\bull) = 
  [\OO_{\Omega_r(E_\bull)}] \cdot
  \Groth_{\wh z(e)}(\breve E_\bull; E_\bull) \,.
\end{equation}
Comparing with (\ref{E:orig_kqf}) it follows that $\Groth_{\wh
  z(r)}(\breve E_\bull; E_\bull) = \big(\sum_\mu c_\mu(r) \prod_i
\Groth_{\mu_i}(E_i-E_{i-1})\big) \cdot \Groth_{\wh z(e)}(\breve
E_\bull; E_\bull)$ on all varieties $X$, so this must hold as an
identity of polynomials in the exterior powers of the bundles $E_i$.
The theorem follows from this.
\ignore{
To deduce the theorem, we use the bundle $Y = \bigoplus_{i=1}^n
\Hom(F_{i-1},F_i)$ over the product of Grassmann bundles
$\Gr(e_0,E_0\oplus \C^M) \times_X \dots \times_X \Gr(e_n,E_n\oplus
\C^M)$ with tautological subbundles $F_0,\dots,F_n$.  The bundle
sequence $E_\bull$ is then the pullback of the universal sequence $F_0
\to \cdots \to F_n$ on $Y$ along a map of varieties $f : X \to Y$.  On
$Y$ we know that $\Groth_{\wh z(r)}(\breve F_\bull; F_\bull) =
[\OO_{\Omega_r(F_\bull)}] \cdot \Groth_{\wh z(e)}(\breve F_\bull;
F_\bull)$.  By taking $M$ sufficiently large, we may therefore assume
that the Grothendieck class of $\Omega_r(F_\bull)$ agrees with
$\Groth_{\wh z(r)}(\breve F_\bull ; F_\bull) / \Groth_{\wh
  z(e)}(\breve F_\bull ; F_\bull)$ in $K(Y)/I$, where $I \subset K(Y)$
is the ideal generated by the classes of subvarieties of codimension
$\dim(X)$ or higher.  The theorem follows from this because
$f^*[\OO_{\Omega_r(F_\bull)}] = [\OO_{\Omega_r(E_\bull)}]$ and $f^*(I)
= 0$.
}
\end{proof}

\begin{remark}
  The above proof cites the main theorem of \cite{buch:grothendieck}
  for the existence of a universal polynomial that expresses the class
  of a quiver variety in the Grothendieck ring of an arbitrary variety
  $X$.  When $X$ has an ample line bundle $L$, one can also deduce the
  theorem directly from (\ref{E:ratio}) as follows.  By twisting the
  sequence $E_\bull$ with a power of $L$, one may assume that all
  bundles $E_i$ are globally generated.  In this case, one can
  construct a bundle $Y = \bigoplus_{i=1}^n \Hom(F_{i-1},F_i)$ over a
  product of Grassmannians $\prod_{i=0}^n \Gr^{e_i}(\C^k)$ with
  tautological quotient bundles $F_0,\dots,F_n$, such that the
  sequence $E_\bull$ is the pullback of the universal sequence $F_0\to
  \dots \to F_n$ on $Y$ along a map of varieties $f: X \to Y$.  On $Y$
  we know that $\Groth_{\wh z(r)}(\breve F_\bull; F_\bull) =
  [\OO_{\Omega_r(F_\bull)}] \cdot \Groth_{\wh z(e)}(\breve F_\bull;
  F_\bull)$.  By taking $k$ sufficiently large, we may therefore
  assume that the Grothendieck class of $\Omega_r(F_\bull)$ agrees
  with $\Groth_{\wh z(r)}(\breve F_\bull ; F_\bull) / \Groth_{\wh
    z(e)}(\breve F_\bull ; F_\bull)$ in $K(Y)/I$, where $I \subset
  K(Y)$ is the ideal generated by the classes of subvarieties of
  codimension $\dim(X)$ or higher.  The theorem follows from this
  because $f^*[\OO_{\Omega_r(F_\bull)}] = [\OO_{\Omega_r(E_\bull)}]$
  and $f^*(I) = 0$.
\end{remark}


\section{Restricted FK-graphs}
\label{sec:restrict}

We will say that an FK-graph $D \subset \N \times \N$ is {\em
  restricted\/} w.r.t.\ the dimension vector $e = (e_0,\dots,e_n)$ if
for every point $(p,q) \in D$ and $0 \leq i \leq n$ we have $p \leq
r_{i-1,0}$ or $q \leq r_{n,i+1}$.  

For each $0 \leq i \leq n$ we let $a^i = (a^i_1,\dots,a^i_{e_i})$ be a
set of $e_i$ variables.  Set $a = (a^0, a^1, \dots, a^n) =
(a^0_1,\dots,a^0_{e_0}, \dots, a^n_1,\dots,a^n_{e_n})$ and $\breve a =
(a^n, a^{n-1}, \dots, a^0) = (a^n_1,\dots,a^n_{e_n}, \dots,
a^0_1,\dots,a^0_{e_0})$.  We need the following variation of
\refcor{C:fkgraph}.

\begin{cor} \label{C:restrict}
For $w \in S_N$ we have
\[ \Groth_w(\breve a; a) =
   \sum_D \,(-1)^{|D|-\ell(w)} \prod_{(p,q) \in D} 
   \left( 1 - \frac{a_p}{\breve a_q} \right)
\]
where the sum is over all FK-graphs $D$ for $w$, which are restricted
{w.r.t.\ }$e$.
\end{cor}
\begin{proof}
  By \refthm{T:stddiag}, the Grothendieck polynomial
  $\Groth_w(\breve a; a)$ is the coefficient of $w$ in the FK-product
  of the diagram $\D'_N$ obtained from $\D_N$ by replacing the top
  variables with $\breve a$ and the left side variables with $a$.  For
  each crossing position $(p,q)$ in this diagram such that $a_p =
  \breve a_q$, the corresponding factor $(1 + (1-\frac{a_p}{\breve
    a_q}) s_{p+q-1})$ of $\Groth(\D'_N)$ is equal to one.  Therefore
  $\Groth(\D'_N)$ is equal to the FK-product of the diagram obtained
  from $\D'_N$ by replacing these crossings with ``\,$\avoid$\,''
  symbols.  Now this diagram is the first of the following two
  equivalent diagrams:
\[ \raisebox{-60pt}{\picT{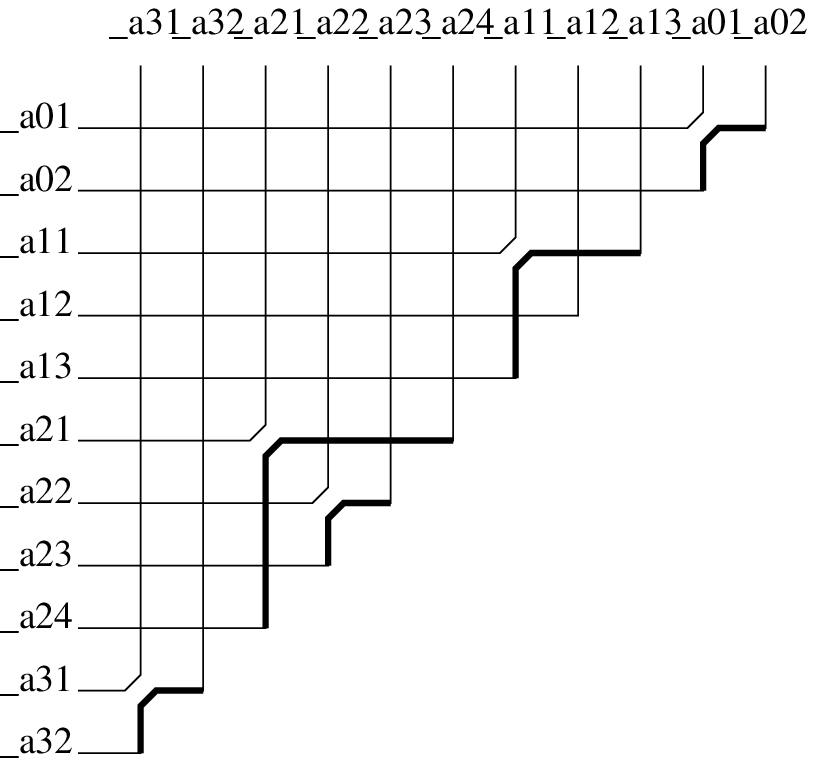}} 
   \approx \ \ \ 
   \raisebox{-60pt}{\picT{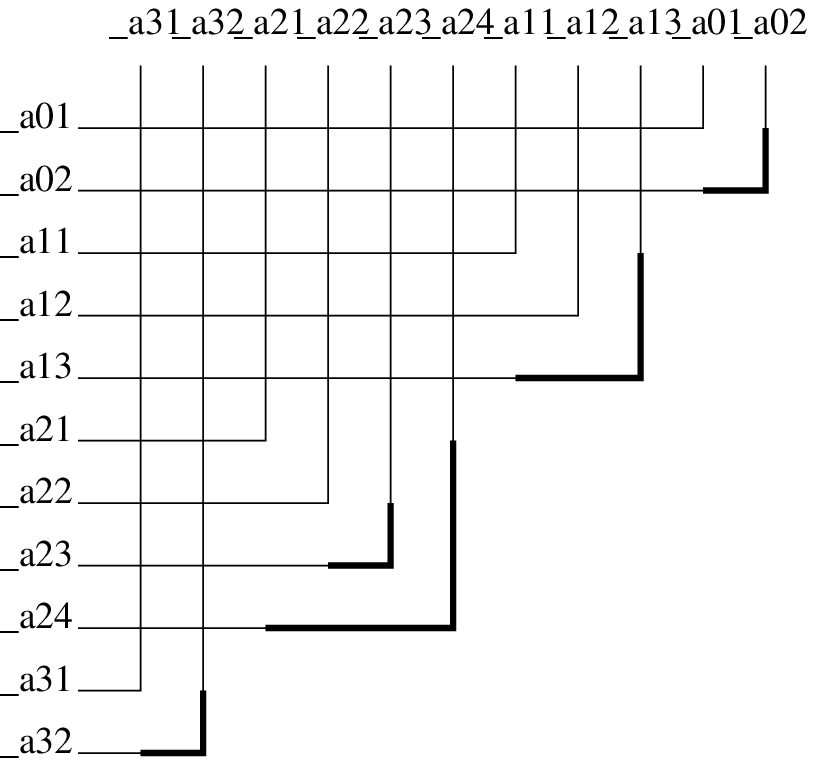}} 
\]
The equivalence follows by using the ``$x=y$'' case of the second
transformation of (\ref{E:key}) to move the thickened line segments
in south-east direction.  The corollary follows from eqn.\ 
(\ref{E:fkprd}) because $\Groth_w(\breve a;a)$ is the coefficient of
$w$ in the FK-product of the second diagram.
\end{proof}


\section{A $K$-theoretic component formula}
\label{sec:component}

For each $1 \leq j \leq n-1$ we set $\delta_j = W_{jj} W_{j+1,j}
\cdots W_{n-1,j} \in S_{r_{n,j-1}}$.  This is the Grassmannian
permutation given by $\delta_j(p) = p + e_{j-1}$ for $e_j < p \leq
r_{nj}$ and $\delta_j(p) = p - r_{n,j+1}$ for $r_{nj} < p \leq
r_{n,j-1}$.  We define a {\em KMS-factorization\/} for the rank
conditions $r$ to be a sequence $(w_1,\dots,w_n)$ of permutations with
$w_i \in S_{e_{i-1}+e_i}$, such that the conjugate Zelevinsky
permutation $z(r)$ is equal to the absolute Hecke product
\begin{equation} \label{E:kmsfac}
  w_1 \cdot \delta_1 \cdot w_2 \cdot \delta_2 
  \cdots \delta_{n-1} \cdot w_n \,.
\end{equation}
In the reduced case, these factorizations are equivalent to the
minimal lace diagrams of Knutson, Miller, and Shimozono
\cite{knutson.miller.ea:four}.\footnote{All permutations $w_i$ are
  inverted in the notation of \cite{knutson.miller.ea:four}.}

A {\em partial permutation\/} contained in the rectangle $k \times l$
with $k$ rows and $l$ columns is a permutation $u \in S_{k+l}$, such
that all descent positions of $u$ are less than or equal to $l$, while
all descent positions of $u^{-1}$ are less than or equal to $k$.  If
this is true, then all FK-graphs for $u$ will be be contained in
$[1,k] \times [1,l]$.

\begin{lemma} \label{L:partial}
  If $(w_1,\dots,w_n)$ is a KMS-factorization for the rank conditions
  $r$, then each permutation $w_j$ is a partial permutation contained
  in $e_{j-1} \times e_j$.
\end{lemma}
\begin{proof}
  Since the absolute Hecke product $\alpha = w_{j+1} \cdot
  \delta_{j+1} \cdot w_{j+2} \cdots \delta_{n-1} \cdot w_n$ is a
  permutation in $S_{r_{nj}}$, we have $\delta_j \cdot \alpha(p) =
  \delta_j(p) = p - r_{n,j+1}$ for $r_{nj} < p \leq r_{n,j-1}$.  If
  $w_j$ had a descent in the interval $[e_j+1,e_j+e_{j-1}-1]$ then the
  product (\ref{E:kmsfac}) would obtain a descent in the interval
  $[r_{nj}+1, r_{n,j-1}-1]$, a contradiction.  By using that
  $(w_n^{-1},\dots,w_1^{-1})$ is a KMS-factorization for the mirrored
  rank conditions $r'_{ij} = r_{n-j,n-i}$, we similarly deduce that
  $w_j^{-1}$ has no descent positions in the interval
  $[e_{j-1}+1,e_{j-1}+e_j-1]$.
\end{proof}

If $u$ is a partial permutation, and if the rectangle $k \times l$ is
understood, we set $\wh u = w_0^{(k+l)} u^{-1} w_0^{(k+l)}$, where
$w_0^{(k+l)}$ is the longest element in $S_{k+l}$.  Notice that the
$180^\circ$ rotation of an FK-graph $D$ for $u$ will be an FK-graph
for $\wh u$.  We will denote this rotated FK-graph by $\wh D$, that
is, $\wh D = \{ (k+1-p,l+1-q) \mid (p,q) \in D \}$.

Given a sequence $(P_1,\dots,P_n)$ of FK-graphs such that each $P_i$
is contained in $[1,e_{i-1}] \times [1,e_i]$, we let
$\wh\Phi(P_1,\dots,P_n)$ denote the FK-graph which is the union of the
sets $\{ (p+r_{i-2,0}, q+r_{n,i+1}) \mid (p,q) \in \wh{P_i} \}$, for
all $1 \leq i \leq n$, together with the unique FK-graph $D_e$ for the
minimal Zelevinsky permutation $\wh z(e)$.
\[ \wh\Phi(P_1,\dots,P_n) \ \ \ = \ \ \ \ \raisebox{-55pt}{\pic{50}{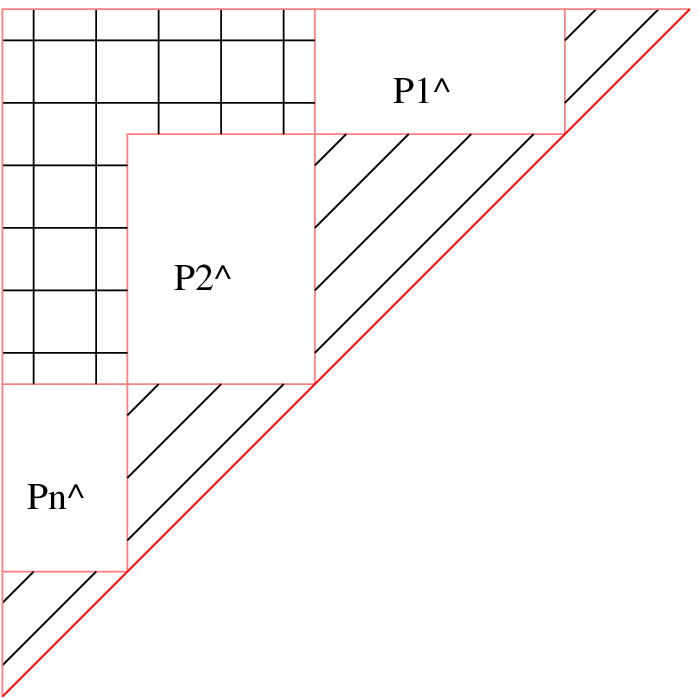}}\]
This construction was used in \cite{knutson.miller.ea:four} (in the
reverse direction) and in \cite{yong:embedding} (for certain special
RC-graphs).  Notice that the crossing positions to the left of
$\wh{P_i}$ in $\wh\Phi(P_1,\dots,P_n)$ form a reduced FK-graph for
$w_0^{(N)} \delta_i^{-1} w_0^{(N)}$, and $D_e$ is the union of these
positions.

\begin{lemma} \label{L:biject}
  Let $r$ be a set of rank conditions.  The map $\wh\Phi$ gives a
  bijection of the set of sequences $(P_1,\dots,P_n)$ of FK-graphs for
  which $(w(P_1),\dots,w(P_n))$ is a KMS-factorization for $r$, with
  the set of all restricted FK-graphs for $\wh z(r)$.
\end{lemma}
\begin{proof}
  The south-west to north-east absolute Hecke product of the simple
  reflections $s_{p+q-1}$ for $(p,q)$ in the block of $\wh{P_i}$ is
  equal to $1^{N-e_i-e_{i-1}} \times w(\wh P_i) = w_0^{(N)}
  w(P_i)^{-1} w_0^{(N)}$.  It follows from this that
\[ w_0^{(N)} w(\wh \Phi(P_1,\dots,P_n))^{-1} w_0^{(N)}
   = w(P_1) \cdot \delta_1 \cdot w(P_2) \cdot \delta_2 \cdots
   \delta_{n-1} \cdot w(P_n) \,.
\]
The lemma therefore follows from \reflemma{L:partial} together with
the definition (\ref{E:kmsfac}) of KMS-factorizations.
\end{proof}

We can now prove the following $K$-theoretic generalization of
\cite[Cor.~6.15]{knutson.miller.ea:four}.

\begin{thm} \label{T:component}
For any set of rank conditions $r$ we have
\[ \frac{\Groth_{\wh z(r)}(\breve a; a)}{\Groth_{\wh z(e)}(\breve a; a)} =
   \sum_{(w_1,\dots,w_n)} (-1)^{\ell(w_1 \cdots w_n)-d(r)} \,
   \Groth_{w_1}(a^1;a^0) \, \Groth_{w_2}(a^2;a^1) \cdots
   \Groth_{w_n}(a^n;a^{n-1})
\]
where the sum is over all KMS-factorizations $(w_1,\dots,w_n)$ for
$r$.
\end{thm}
\begin{proof}
  \refcor{C:restrict} and \reflemma{L:biject} imply that ${\Groth_{\wh
      z(r)}(\breve a; a)}/{\Groth_{\wh z(e)}(\breve a; a)}$ is equal
  to
\begin{multline*} \sum (-1)^{\ell(w_1 \cdots w_n)-d(r)} \,
   \Groth_{w_1}(a^1_{e_1},\dots,a^1_1; a^0_{e_0},\dots,a^0_1) \cdot \\
   \Groth_{w_2}(a^2_{e_2},\dots,a^2_1; a^1_{e_1},\dots,a^1_1) \cdots
   \Groth_{w_n}(a^n_{e_n},\dots,a^n_1;a^{n-1}_{e_{n-1}},\dots,a^{n-1}_1)
   \,.
\end{multline*}
The theorem follows because $\Groth_{\wh z(r)}(\breve a; a)$ is
symmetric in each set of variables $a^i$.
\end{proof}

\begin{remark}
  Let $(w_1,\dots,w_n)$ be a sequence of partial permutations with
  $w_i$ contained in $e_{i-1} \times e_i$.  Suppose $1 \leq j < n$ and
  $1 \leq k < e_j$ are given such that $w_j(k) < w_j(k+1)$ and
  $w_{j+1}^{-1}(k) < w_{j+1}^{-1}(k+1)$.  If any of the sequences
  $(w_1, \dots, w_j s_k, w_{j+1}, \dots, w_n)$, $(w_1, \dots, w_j, s_k
  w_{j+1}, \dots, w_n)$, or $(w_1, \dots, w_j s_k, s_k w_{j+1}, \dots,
  w_n)$ is a KMS-fac\-tori\-zation, then the definition shows that all
  three are KMS-factorizations.  These transformations were first
  observed during an attempt to generalize the symmetry arguments of
  \cite{buch.feher.ea:positivity} to $K$-theory.  In fact, in
  \cite{buch.feher.ea:positivity} it is proved that all
  KMS-factorizations of a given Zelevinsky permutation are connected
  by these transformations, which gives an easy way to find all of
  them.
\end{remark}


\section{KMS-factorizations are factor sequences}
\label{sec:facseq}

We define the {\em permutation diagram\/} for the rank conditions $r$
to be the part of the matrix (\ref{E:permmat}) which is on or below
the antidiagonal.
\newcommand{\adots}{\raisebox{1pt}{.}\hspace{1pt}\raisebox{4pt}{.}\hspace{1pt}\raisebox{7pt}{.}}
\begin{equation} \label{E:permdiag} \begin{matrix}
         &        &        & W_{n-1,n} \\
         &        & \adots & \vdots \\
         & W_{12} & \dots  & W_{1n} \\
  W_{01} & W_{02} & \dots  & W_{0n} 
\end{matrix} \end{equation}
This diagram can also be obtained from the rectangle diagram of
\cite[\S 2.1]{buch.fulton:chern} by replacing each rectangle $R_{ij}$
with the corresponding Grassmannian permutation $W_{ij}$ (and rotating
the result 45 degrees counter clockwise.)

For $0 \leq k \leq n$ we let $r^{(k)}$ denote the rank conditions
obtained from $r$ by dropping all integers $r_{ij}$ with $j-i < k$,
that is $r^{(k)} = \{ r^{(k)}_{ij} \mid 0 \leq i \leq j \leq n-k \}$
where $r^{(k)}_{ij} = r_{i,j+k}$.  The permutation diagram for these
rank conditions is obtained from (\ref{E:permdiag}) by dropping the
top $k$ diagonals of matrices $W_{ij}$ with $j-i \leq k$.

In this section we prove that all KMS-factorizations for $r$ can be
obtained from the permutation diagram in the following way.  First, if
$n=1$ then the diagram has only one permutation $W_{0n}$, and the only
KMS-factorization is $(W_{0n})$.  Otherwise any KMS-factorization can
be obtained by first constructing a KMS-factorization
$(\alpha_1,\dots,\alpha_{n-1})$ for $r^{(1)}$ and choosing arbitrary
factorizations $\alpha_i = u_i \cdot v_i$ w.r.t.\ the absolute Hecke
product.  Then the sequence $(W_{01}\cdot u_1, v_1\cdot W_{12} \cdot
u_2, \dots, v_{n-1}\cdot W_{n-1,n})$ is a KMS-factorization
for the rank conditions $r$.  The exact statement that we prove is the
following theorem, which also includes a criterion for
KMS-factorizations similar to the criterion for factor sequences
proved in \cite{buch:on}.

\begin{thm} \label{T:facseq}
  (a) If $(w_1,\dots,w_n)$ is a KMS-factorization for $r$, then each
  permutation $w_i$ has a reduced factorization $w_i = v_{i-1} \cdot
  W_{i-1,i} \cdot u_i$ with $v_{i-1} \in S_{e_{i-1}}$ and $u_i \in
  S_{e_i}$, such that $v_0 = u_n = 1$.
  
  (b) Let $u_1,v_1,\dots,u_{n-1},v_{n-1}$ be permutations with
  $u_i,v_i \in S_{e_i}$.  Then the sequence $(W_{01}\cdot u_1, v_1
  \cdot W_{12} \cdot u_2, \dots, v_{n-1} \cdot W_{n-1,n})$ is a
  KMS-factorization for $r$ if and only if $(u_1 \cdot v_1, u_2 \cdot
  v_2, \dots, u_{n-1} \cdot v_{n-1})$ is a KMS-factorization for
  $r^{(1)}$.
\end{thm}

It is also possible to formulate a tableau-based version of this
theorem.  For cohomological quiver coefficients, this is based on
Fomin and Greene's formula for Stanley coefficients
\cite{fomin.greene:noncommutative}, and uses a diagram of tableaux
obtained from (\ref{E:permdiag}) by replacing $W_{ij}$ with the unique
row and column {\em decreasing\/} tableau $T_{ij}$ of shape $R_{ij}$
(i.e.\ with $r_{i+1,j}-r_{ij}$ rows and $r_{i,j-1}-r_{ij}$ columns),
such that the reading word of $T_{ij}$ is a reduced word for $W_{ij}$.
Factor sequences are generated using the Coxeter-Knuth product, and it
follows from \refthm{T:facseq}, \refcor{C:quivcoef}, and
\cite[Thm.~1.2]{fomin.greene:noncommutative} that a cohomological
quiver coefficient $c_\mu(r)$ counts the number of factor sequences of
shape $\mu$.  In the general case there remains work to be done.  We
refer to \cite{bkty:stable} for details.

Let $0 \leq k < n$ and consider permutations $w_1,\dots,w_{n-k}$ such
that $w_i \in S_{r_{i,k+i-1}}$ for all $i$.  We let
$\Phi_k(w_1,\dots,w_{n-k}) \in S_N$ denote the south-west to
north-east product of the matrix obtained from (\ref{E:permmat}) by
replacing $W_{i-1,k+i}$ with $w_i$ for $1 \leq i \leq n-k$, and by
replacing $W_{ij}$ with the identity for $j-i \geq k+2$.  When $n=4$
and $k=1$ this matrix looks as follows:
\[ \begin{matrix}
        W_{31} & W_{32} & W_{33} & W_{34} \\
        W_{21} & W_{22} & W_{23} &  w_3   \\
        W_{11} & W_{12} &  w_2   &   1    \\
        W_{01} &  w_1   &   1    &   1
\end{matrix} \]
Notice that $(w_1,\dots,w_n)$ is a KMS-factorization for $r$ if and
only if $\Phi_0(w_1,\dots,w_n) = z(r)$.  Furthermore, if $k > 0$ and
$u_1,v_1,\dots,u_{n-k},v_{n-k}$ are permutations such that $u_i,v_i
\in S_{r_{i,k+i-1}}$ then
\[ \Phi_k(u_1 \cdot v_1, \dots, u_{n-k} \cdot v_{n-k}) =
   \Phi_{k-1}(W_{0,k} \cdot u_1, v_1 \cdot W_{1,k+1} \cdot u_2,
   \dots, v_{n-k} \cdot W_{n-k,n}) \,.
\]

\begin{lemma} \label{L:factor}
  If $\Phi_k(w_1,\dots,w_{n-k}) = z(r)$ then each permutation $w_i$
  has a reduced factorization $w_i = v \cdot W_{i-1,k+i} \cdot u$,
  where $v \in S_{r_{i-1,k+i-1}}$ and $u \in S_{r_{i,k+i}}$.
  Furthermore, $v$ is trivial when $i=1$ while $u$ is trivial when
  $i=n-k$.
\end{lemma}
\begin{proof} Fix $i$ and set $j = k+i-1$, $a = r_{i-1,j}$, $b =
  r_{i,j+1}$, and $c = r_{i-1,j+1}$.  For the first assertion it is
  enough to show that $w_i \in S_{a+b-c}$, that $w_i(p) \leq a$ for $b
  < p \leq a+b-c$, and that $w_i^{-1}(p) \leq b$ for $a < p \leq
  a+b-c$.  In fact, $v w_i u$ will then satisfy the same properties
  for every $v \in S_a$ and $u \in S_b$.  If we choose $v$ and $u$
  such that $v^{-1} w_i u^{-1}$ has no descents in the interval
  $[1,b-1]$, $u w_i^{-1} v$ has no descents in the interval $[1,a-1]$,
  and $\ell(v^{-1} w_i u^{-1}) = \ell(w_i) - \ell(v) - \ell(u)$, then
  we must have $v^{-1} w_i u^{-1} = W_{i-1,j+1}$, so we can use the
  factorization $w_i = v \cdot W_{i-1,j+1} \cdot u$.  
  
  Let $\sigma = W_{i-1,1} W_{i-1,2} \cdots W_{i-1,j}$ be the
  product of the permutations west of $w_i$ in the matrix
  (\ref{E:permmat}), and $\tau = \prod_{s=i}^{n-1} \prod_{t=1}^{j+1}
  W_{st}$ the product of the permutations weakly west and strictly
  north of $w_i$ in this matrix.  It follows by induction on $j$ that
  all descent positions of $\sigma$ are greater than or equal to
  $r_{ij}$.  A similar argument shows that the descent positions of
  $\tau^{-1}$ are greater than or equal to $r_{ij}$.  Furthermore, the
  product $\sigma \tau$ is reduced, since it is part of the defining
  reduced factorization of the conjugate Zelevinsky permutation for
  rank conditions obtained by replacing $r_{i-1,j+1}$ with
  $\min(r_{i-1,j}, r_{i,j+1})$.  Since $w_i \in S_{r_{ij}}$, it
  follows that $\sigma w_i \tau$ is also a reduced product.  We can
  therefore write $z(r) = \alpha \sigma w_i \tau \beta$ as a product
  of permutations, where $\alpha \in S_{r_{i-1,0}}$ and $\beta \in
  S_{r_{n,j+1}}$.
  
  Notice that for $b < p \leq r_{ij}$ we have $\tau
  \beta(r_{n,j+1}-b+p) = \tau(r_{n,j+1}-b+p) = p$, and for $a < p \leq
  r_{ij}$ we have $\alpha \sigma(p) = \alpha(r_{i-1,0}-a+p) =
  r_{i-1,0}-a+p$.  For $a+b-c < p \leq r_{ij}$ it follows from
  \reflemma{L:zeldesc} that $z(r)(r_{n,j+1}-b+p) = r_{i-1,0}-a+p$, so
  we must have $w_i(p) = p$, that is $w_i \in S_{a+b-c}$.  The lemma
  also shows that $z(r)(r_{n,j+1}-b+p) \leq r_{i-1,0}$ for all $b < p
  \leq a+b-c$.  This implies that $w_i(p) \leq a$, since otherwise we
  would have $\alpha \sigma w_i \tau \beta(r_{n,j+1}-b+p) = \alpha
  \sigma w_i(p) = r_{i-1,0}-a+w_i(p) > r_{i-1,0}$.  A symmetric
  argument shows that $w_i^{-1}(p) \leq b$ for $a < p \leq a+b-c$,
  which completes the proof of the first assertion.
  
  When $i=n-k$, $\beta$ is trivial, so we have $\tau \beta(p) = p$ for
  all $p \leq b = r_{in}$.  If $u \in S_b$ was not the identity then
  $z(r)$ would get a descent in the interval $[1,e_n-1]$, a
  contradiction.  A similar argument shows that $v$ must be trivial
  when $i=1$.
\end{proof}

\extra{Suppose that all descent positions of $\sigma$ and $\tau^{-1}$
  are greater than or equal to $m$ and that $\sigma \tau$ is a reduced
  product.  We claim that for any $w \in S_m$, the product $\sigma w
  \tau$ is also reduced.  Since $\sigma w$ is clearly reduced, this
  means that we can find $p < q$ such that $\sigma w(p) > \sigma w(q)$
  and $\tau^{-1}(p) > \tau^{-1}(q)$.  Since $\sigma \tau$ is reduced
  we must have $p \leq m$, and the assumption about the descents of
  $\tau^{-1}$ implies that $q > m$.  Set $p' = \max(p,w(p))$.  We then
  have $p' \leq m < q$, $\sigma(p') \geq \sigma(w(p)) > \sigma(q)$,
  and $\tau^{-1}(p') \geq \tau^{-1}(p) > \tau^{-1}(q)$, so $\sigma
  \tau$ is not a reduced product, a contradiction.}

\begin{proof}[Proof of Theorem \ref{T:facseq}]
  Part (a) is a special case of \reflemma{L:factor} so we prove part
  (b).  By induction on $n$ we may assume that the theorem is true for
  all rank conditions $r^{(k)}$ with $k \geq 1$.  We claim that
  $(w_1,\dots,w_{n-k})$ is a KMS-factorization for $r^{(k)}$ if and
  only if $\Phi_k(w_1,\dots,w_{n-k}) = z(r)$.  If either is true, then
  we can write $w_i = v_{i-1} \cdot W_{i-1,k+i} \cdot u_i$ for each
  $i$, where $u_i, v_i \in S_{r_{i,k+i}}$, and $v_0 = u_{n-k} = 1$.
  This proves the claim when $k=n-1$.  For $1 \leq k \leq n-2$ we know
  by induction on $k$ that $\Phi_k(w_1,\dots,w_{n-k}) = \Phi_{k+1}(u_1
  \cdot v_1, \dots, u_{n-k-1} \cdot v_{n-k-1})$ equals $z(r)$ if and
  only if $(u_1 \cdot v_1, \dots, u_{n-k-1} \cdot v_{n-k-1})$ is a
  KMS-factorization for $r^{(k+1)}$, which by the theorem for
  $r^{(k)}$ is equivalent to $(w_1,\dots,w_{n-k})$ being a
  KMS-factorization for $r^{(k)}$.  This proves the claim, and the
  theorem follows from the claim with $k=1$.
\end{proof}

In the next section we need the following generalization of
\cite[Cor.~4.12]{knutson.miller.ea:four}.

\begin{cor} \label{C:kms_stab}
  The KMS-factorizations for $r+m$ are precisely the sequences\linebreak
  $(1^m\times w_1, \dots, 1^m \times w_n)$ for which $(w_1,\dots,w_n)$
  is a KMS-factorization for $r$.
\end{cor}
\begin{proof} 
  This is immediate from \refthm{T:facseq} because the
  permutation diagram for $r+m$ is obtained from (\ref{E:permdiag}) by
  replacing each permutation $W_{ij}$ with $1^m \times W_{ij}$.
\end{proof}


\section{Alternating signs of quiver coefficients}
\label{sec:altsigns}

We can now prove the stable version of the component formula.  We note
that \refthm{T:rankstab} of the next section can be substituted for
the reference to geometry in its proof.

\begin{thm} \label{T:stable}
  For any set of rank conditions $r$ we have
\[ \frac{\Groth_{\wh z(r)}(\breve a; a)}{\Groth_{\wh z(e)}(\breve a; a)} =
   \sum_{(w_1,\dots,w_n)} (-1)^{\ell(w_1 \cdots w_n)-d(r)} \,
   G_{w_1}(a^1;a^0) \, G_{w_2}(a^2; a^1) \cdots G_{w_n}(a^n;a^{n-1})
\]
where the sum is over all KMS-factorizations for $r$.
\end{thm}
\begin{proof}
  It follows from (\ref{E:orig_kqf}) and \refthm{T:ratio} that
  ${\Groth_{\wh z(r)}(\breve a; a)}/{\Groth_{\wh z(e)}(\breve a; a)}$
  is equal to $\sum_\mu c_\mu(r)\, G_{\mu_1}(a^1;a^0) \cdots
  G_{\mu_n}(a^n;a^{n-1})$.  Since $c_\mu(r) = c_\mu(r+m)$ we deduce
  that
\[ \frac{\Groth_{\wh z(r)}(\breve a; a)}{\Groth_{\wh z(e)}(\breve a; a)} 
   = \frac{\Groth_{\wh z(r+m)}(a^n,1^m, \dots, a^0,1^m \ ;\  
   a^0,1^m, \dots, a^n,1^m)}{\Groth_{\wh z(e+m)}( a^n,1^m, \dots, a^0,1^m
   \ ;\ a^0,1^m, \dots, a^n,1^m )}
\]
for all integers $m \geq 0$.  For $m \geq \max(e_0,\dots,e_n)$ it
follows from \refthm{T:component} and \refcor{C:kms_stab} that the
right hand side of this identity equals
\[\begin{split}
  & \sum (-1)^{\ell(w_1\dots w_n)-d(r)}\, \Groth_{1^m \times w_1}(a^1,1^m \,;\,
  a^0,1^m) \cdots \Groth_{1^m \times w_n}(a^n,1^m \,;\, a^{n-1},1^m) \\
  & = \ \sum (-1)^{\ell(w_1\dots w_n)-d(r)}\, G_{w_1}(a^1;a^0) \cdots 
  G_{w_n}(a^n;a^{n-1})
\end{split}\]
as required.
\ignore{
  Fix an integer $m \geq \max(e_0,\dots,e_n)$.  By using the geometric
  fact that $c_\mu(r) = c_\mu(r+m)$, as well as \refthm{T:component}
  and \refcor{C:kms_stab}, we then obtain
\[ \begin{split} 
& \frac{\Groth_{z(r)}(\breve a; a)}{\Groth_{z(e)}(\breve a; a)} 
  = \frac{\Groth_{z(r+m)}(a^n,1^m, \dots, a^0,1^m \ ;\  
    a^0,1^m, \dots, a^n,1^m)}{\Groth_{z(e+m)}( a^n,1^m, \dots, a^0,1^m \ ;\  
    a^0,1^m, \dots, a^n,1^m )} \\
&\ \ \ = 
  \sum (-1)^{\ell(w_1\dots w_n)-d(r)}\, \Groth_{1^m \times w_1}(a^1,1^m \,;\,
  a^0,1^m) \cdots \Groth_{1^m \times w_n}(a^n,1^m \,;\, a^{n-1},1^m) \\
&\ \ \ = 
  \sum (-1)^{\ell(w_1\dots w_n)-d(r)}\, G_{w_1}(a^1;a^0) \cdots 
  G_{w_n}(a^n;a^{n-1})
\end{split} \]
as required.
}
\end{proof}

It was conjectured in \cite{buch:grothendieck} that quiver
coefficients have signs which alternate with codimension.  This
conjecture is a consequence of the following explicit formula for
quiver coefficients, which follows from \refthm{T:stable} and
\cite[Thm.~4]{lascoux:transition}.

\begin{cor} \label{C:quivcoef}
The quiver coefficient $c_\mu(r)$ is given by
\[ c_\mu(r) = (-1)^{\sum |\mu_i| - d(r)}\, \sum_{(w_1,\dots,w_n)}
   | c_{w_1,\mu_1} c_{w_2,\mu_2} \cdots c_{w_n,\mu_n} |
\]
where the sum is over all KMS-factorizations for the rank conditions
$r$, and $c_{w_i,\mu_i}$ is defined by equation (\ref{E:stabcoef}).
\end{cor}


\section{Double quiver polynomials}
\label{sec:double}

Let $r$ be a set of rank conditions.  For each $0 \leq i \leq n$ we
let $b^i = (b^i_1, \dots, b^i_{e_i})$ be a set of $e_i = r_{ii}$
variables, and we set $b = (b^0,b^1,\dots,b^n)$.  The $K$-theoretic
analogue of the double ratio formula from
\cite{knutson.miller.ea:four} defines the polynomial
\[ K_r(a;b) = \frac{\Groth_{\wh z(r)}(\breve a;b)}
   {\Groth_{\wh z(e)}(\breve a;b)} \,,
\]
which is named a {\em double quiver polynomial\/} in
\cite{knutson.miller.ea:four}.  In this section we prove some facts
about such polynomials.  In particular, we show that $K_r(a;b)$ is a
specialization of the quiver formula constructed in
\cite{buch:grothendieck}.  When the rank conditions have the form
$r+m$ for $m$ large, the cohomology versions of theorems
\ref{T:multisuper}, \ref{T:rankstab}, and \ref{T:double} below follow
from the results in \cite{knutson.miller.ea:four} (see Cor.~6.13 and
Thm.~6.20 in loc.\ cit.)

It follows from (i) and (ii) of \reflemma{L:zelprop} that $K_r(a;b)$
is separately symmetric in each set of variables $a^i$ and $b^i$ (and
that the variables $a^0$ and $b^n$ do not occur).  In addition we
have:

\begin{thm} \label{T:multisuper}
The polynomial $K_r(a;b)$ is multi supersymmetric, that is, if one
sets $a^{i+1}_1 = b^i_1$ in $K_r(a;b)$ then the result is independent
of these variables.
\end{thm}

Since $\Groth_{\wh z(e)}(\breve a; b) = \prod_{(p,q) \in D_e} (1 -
\frac{b_p}{\breve a_q})$, this theorem is an immediate consequence of
the following proposition, applied to $\wh z(r)$ with $k =
r_{i-1,0}+1$ and $l = r_{n,i+2}+1$.

\begin{prop}
  Let $w \in S_N$ be a permutation and let $k, l \in \N$.  Assume that
  $w(i) > k$ for $1 \leq i < l$ and that $w^{-1}(i) > l$ for $1 \leq i
  < k$.  Then the polynomial  
\[ \frac{\Groth_w(a_1,\dots,a_{l-1},b_k,a_{l+1},\dots,a_N\,;\,b_1,\dots,b_N)}{
   \left(\prod_{i=1}^{k-1} (1-\frac{b_i}{b_k})\right) 
   \left(\prod_{i=1}^{l-1} (1-\frac{b_k}{a_i})\right)}
\]
is independent of $b_k$.
\end{prop}
\begin{proof}
  By \refthm{T:stddiag}, $\Groth_w(a_1,\dots,a_{l-1},b_k,
  a_{l+1},\dots,a_N ; b_1,\dots,b_N)$ is the coefficient of $w$ in the
  FK-product of the first of the following three equivalent diagrams:
\[ \raisebox{-40pt}{\picT{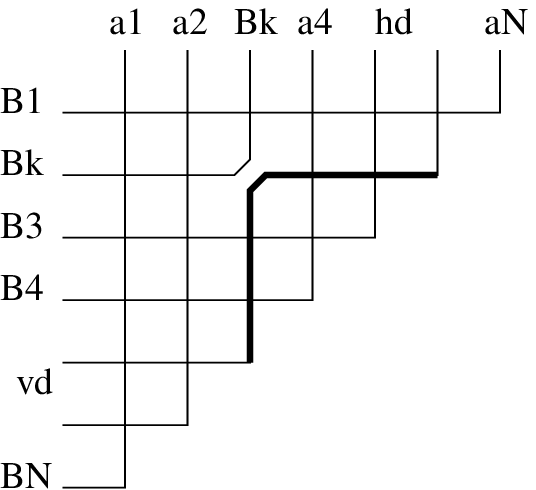}} \approx \ \ 
   \raisebox{-40pt}{\picT{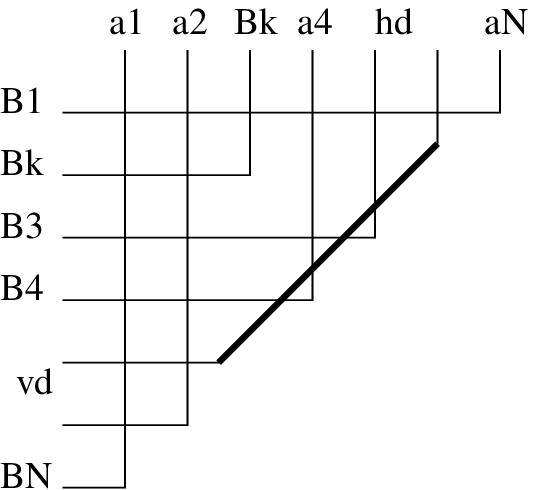}} \approx \ \ 
   \raisebox{-40pt}{\picT{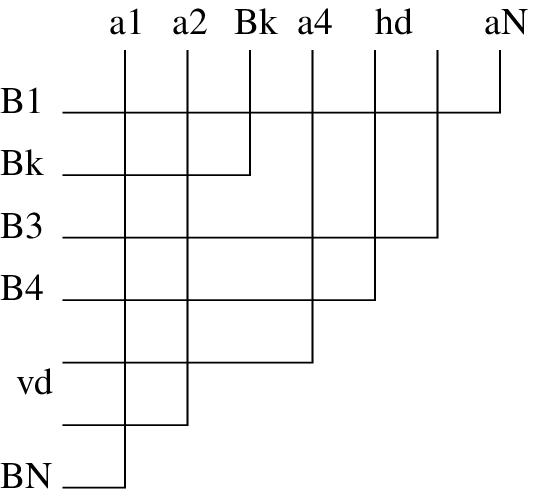}}
\]
The first diagram is obtained from $\D_N$ by replacing $a_l$ with
$b_k$ and the crossing at position $(k,l)$ with ``\,$\avoid$\,'', and
the others are the result of moving the thick line segments labeled
$b_k$ in south-east direction, using the rules of (\ref{E:key}).

Notice that if $D$ is a subset of the crossing positions of the third
diagram such that $w(D) = w$, then $D$ must contain all crossings
involving the line segments labeled $b_k$.  The proposition therefore
follows from eqn.~(\ref{E:fkprd}).
\end{proof}

Our next theorem is a $K$-theoretic generalization of Conjecture~6.14
from \cite{knutson.miller.ea:four}.  It says that the double ratio
formula satisfies rank stability.

\begin{thm} \label{T:rankstab}
  For any rank conditions $r$ and non-negative integer $m$ we have
\[ K_{r+m}(a^0,1^m,\dots,a^n,1^m \,;\, b^0,1^m,\dots,b^n,1^m)
   = K_r(a^0,\dots,a^n \,;\, b^0,\dots,b^n) \,.
\]
\end{thm}
\begin{proof}
Using symmetry, it is enough to prove that
\begin{equation} \label{E:doublestab} \begin{split} 
  K_{r+1}(a^0,1,\, c_1,a^1,\dots,c_n,a^n \,;\,
  c_1,b^0,\dots,c_n,b^{n-1} ,\, b^n,1) \ \ \ \ \ \ \ \ \ \ \ \\
  = \ \  K_r(a^0,\dots,a^n \,;\, b^0,\dots,b^n) \,.
\end{split} \end{equation}
The polynomial $\Groth_{\wh z(r+1)} = \Groth_{\wh z(r+1)}( c_n,a^n,
\dots, c_1,a^1, a^0,1 \,;\, c_1,b^0,\dots,c_n,b^{n-1}, b^n,1)$ is the
coefficient of $\wh z(r+1)$ in the FK-product of the following two
equivalent diagrams (of size $N+n+1$):
\[ \raisebox{-75pt}{\pic{50}{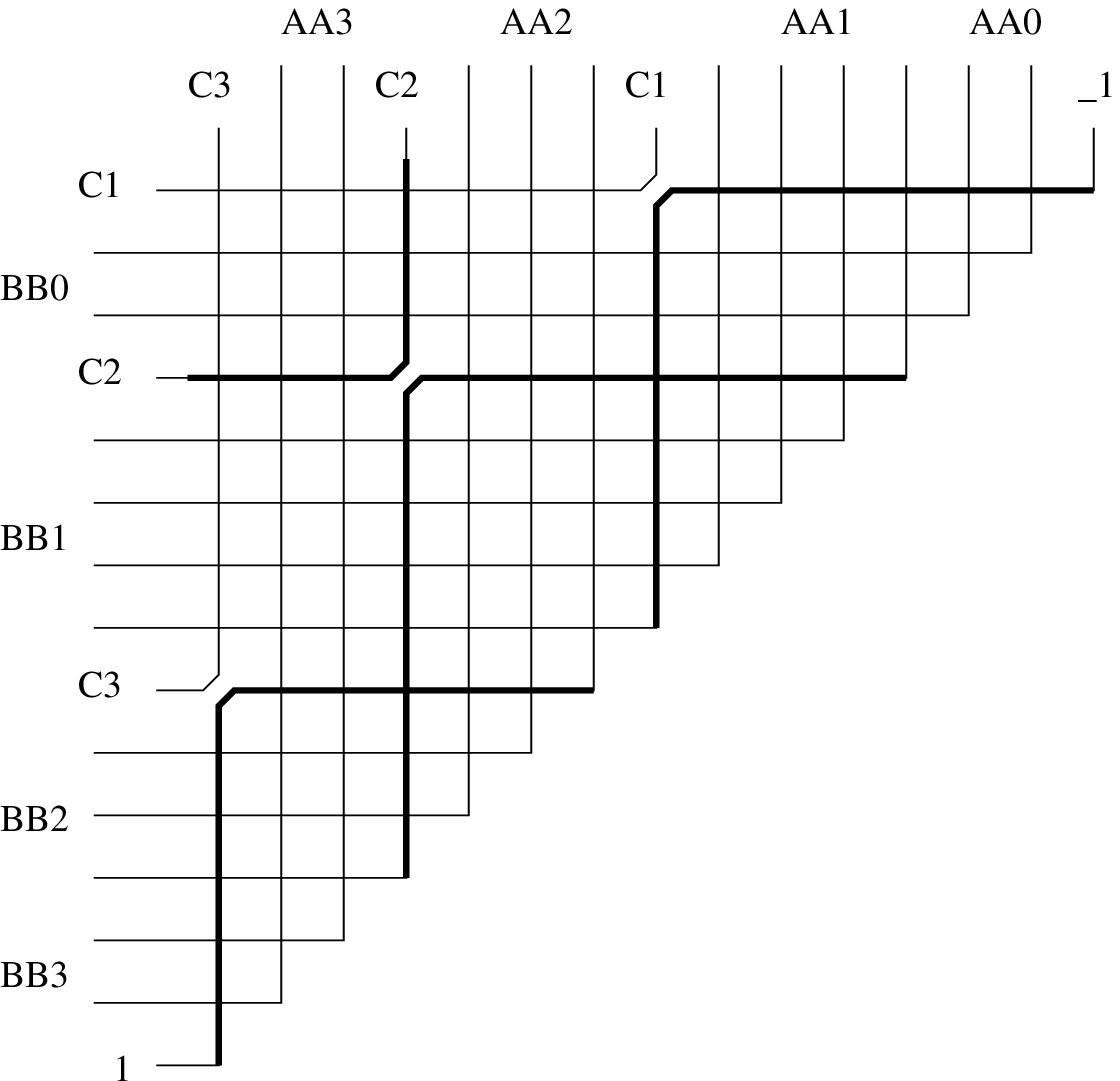}} \approx \ \ 
   \raisebox{-75pt}{\pic{50}{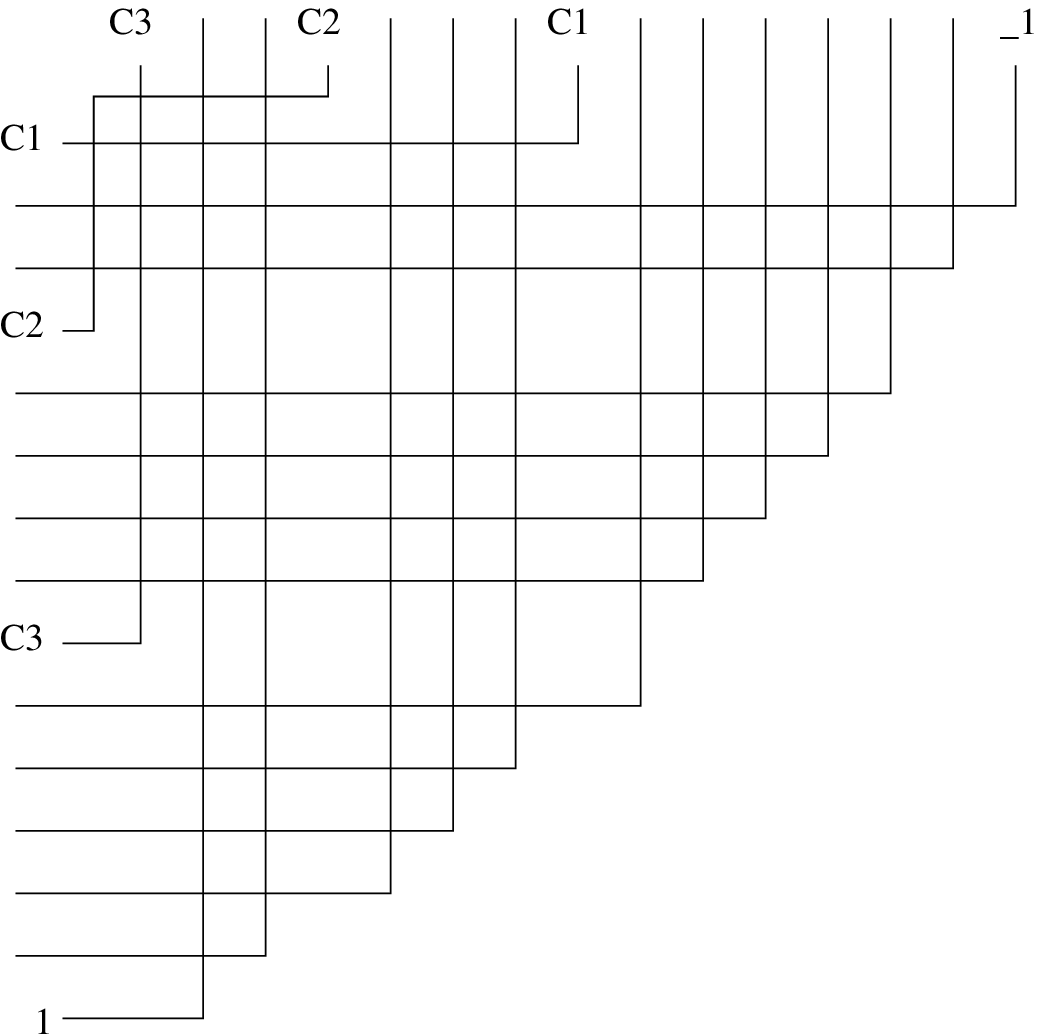}}
\]
As indicated in the first diagram, the line segments extending
furthest north and west are labeled with the variable sets $a^i$ and
$b^i$, respectively.  The equivalence is obtained by moving the
thickened line segments toward the borders of the diagram.  Since $\wh
z(r+1)$ fixes $N+n+1$, it follows that $\Groth_{\wh z(r+1)}$ is also
the coefficient of $\wh z(r+1)$ in the FK-product of the diagram:
\[  \D \ \ = \ \ \ 
  \raisebox{-70pt}{\pic{50}{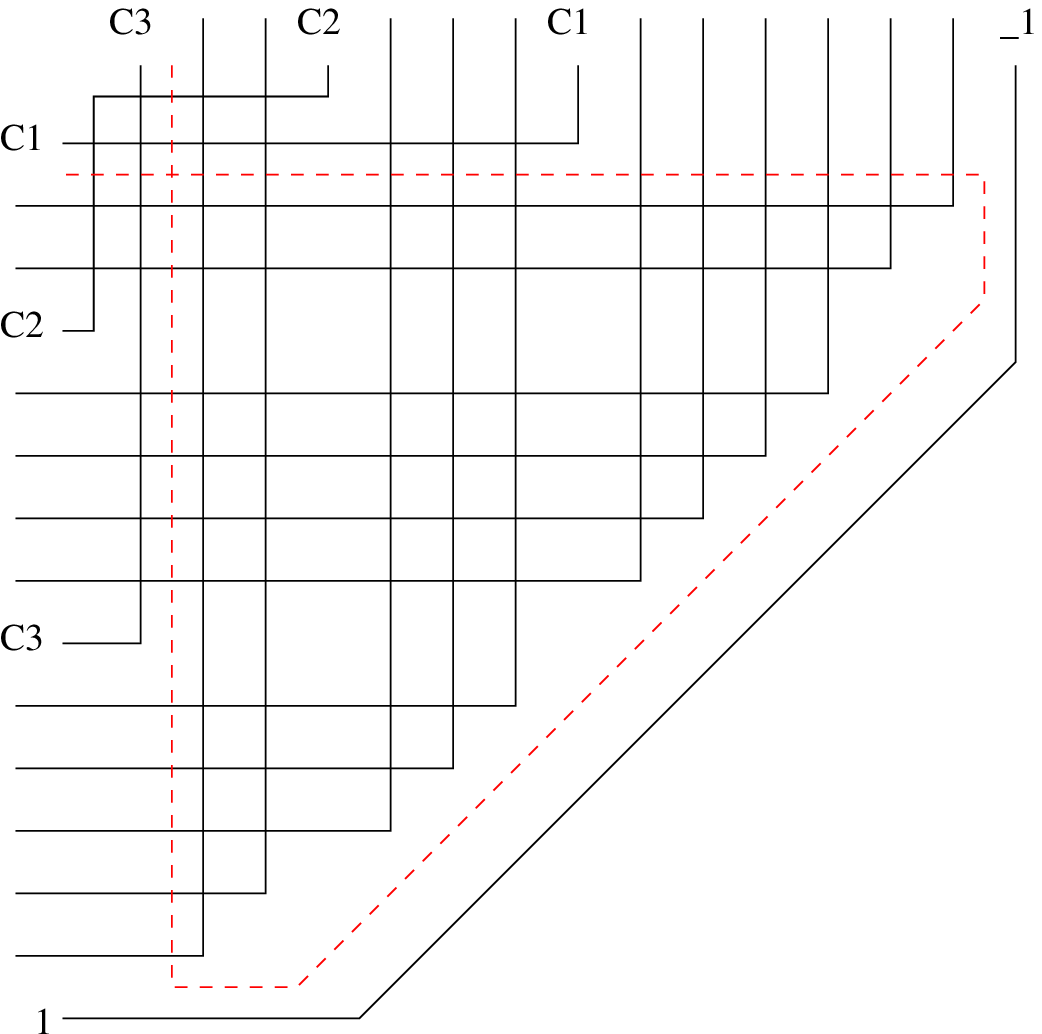}}
\]


We claim that the FK-graphs for $\wh z(r+1)$ w.r.t.\ this diagram are
exactly those obtained by placing an FK-graph for $\wh z(r)$ w.r.t.\ 
$\D_N$ in the triangular region of $\D$, and including all crossing
positions $P$ outside this region.  The theorem follows from this
claim because it shows that $\Groth_{\wh z(r+1)} = {\mathcal Q} \cdot
\Groth_{\wh z(r)}(a^n,a^{n-1},\dots,a^0 \,;\, b^0,b^1,\dots,b^n)$,
where ${\mathcal Q}$ is the product of the factors $h(P)$ for the
crossing positions $P$ outside the triangular region, and
eqn.~(\ref{E:doublestab}) is immediate from this identity.

We will say that a permutation $w \in S_N$ satisfies the condition (*)
for the set of rank conditions $r$, if for all $i$ and $p \leq
r_{n,i+1}$ we have $w(p) > r_{i-1,0}$.  As noted in section
\ref{sec:zel}, $\wh z(r)$ has this property.  Notice that if $D
\subset C(\D)$ is any FK-graph such that $w(D)$ satisfies (*) for
$r+1$, then $D$ must contain all the crossing positions outside the
triangular region.  In particular, all FK-graphs for $\wh z(r+1)$ must
contain these crossing positions.

The crossing positions in the south-west region of $\D$ represent the
Grassmannian permutation $\alpha$ defined by $\alpha(i) = r_{i-2,0}+i$
for $1 \leq i \leq n$ and $\alpha(i) < \alpha(i+1)$ for $i \neq n$,
while the permutation $\beta$ represented by the north-east crossings
is given by $\beta^{-1}(i) = r_{n,n+2-i}+i$ for $1 \leq i \leq n$ and
$\beta^{-1}(i) < \beta^{-1}(i+1)$ for $i \neq n$.  If we let $w \in
S_N$ be the permutation of an FK-graph placed in the triangular region
of $\D$, then the claim states that the absolute Hecke product $\alpha
\cdot (w_0^{(n)} \times w) \cdot \beta$ equals $\wh z(r+1)$ if and
only if $w = \wh z(r)$.

If $\alpha \cdot (w_0^{(n)} \times w) \cdot \beta$ satisfies (*) for
$r+1$, then this product must be reduced.  \mbox{Otherwise} one could
skip a simple reflection factor of $\alpha$ or $\beta$, which could be
exploited to construct an FK-graph for the product that missed the
corresponding crossing outside the triangular region in $\D$.  Since a
similar argument applies if the (usual) product of permutations
$\alpha (w_0^{(n)} \times w) \beta$ satisfies (*) for $r+1$, it is
enough to check that $\alpha (w_0^{(n)} \times w) \beta$ equals $\wh
z(r+1)$ if and only if $w = \wh z(r)$.

Equivalently, we must check that the {\em inverse\/} of $\alpha
(w_0^{(n)} \times \wh z(r)) \beta$ maps $N+n+2-p$ to $N+n+2-z(r+1)(p)$
for every $2 \leq p \leq N+n+1$.  We use \reflemma{L:zeldesc} to do
this.  Set $r^+ = r+1$.  This means that $r^+_{ij} = r_{ij} + 1$ for
$i \leq j$ while $r^+_{ij} = r_{ij} + 1+i-j$ when $i > j$.  Now choose
$0 \leq i \leq j+1 \leq n+1$ such that $r^+_{n,j+1} + r^+_{i-1,j} -
r^+_{i-1,j+1} < p \leq r^+_{n,j+1}+r^+_{i,j}-r^+_{i,j+1}$.

If $0 \leq j \leq n-1$ and $p = r^+_{n,j}$ , then $i = j+1$ and
$z(r^+)(p) = r^+_{j+1,0}$.  Furthermore, $N+n+2-p = r_{j-1,0}+j+1$ is
mapped to $j+1$ by $\alpha^{-1}$, which in turn is mapped to $n-j$ by
$w_0^{(n)}$, and $\beta^{-1}(n-j) = r_{n,j+2}+n-j = N+n+2-r^+_{j+1,0}$
as required.

Otherwise we have $p \neq r^+_{n,j}$ or $j = n$.  Set $i' = \max(i,1)$
and $j' = \min(j,n-1)$.  Then $z(r^+)(p) = p + \delta + i' + j' - n$,
where $\delta = - r_{n,j+1} + r_{i,j+1} +r_{i-1,0} -r_{i-1,j}$.  Since
$j' \leq n-1$ is maximal such that $\alpha(j'+1) < N+n+2-p$, we get
$\alpha^{-1}(N+n+2-p) = (N+n+2-p)+n-j'-1$.  One now checks that
$r_{n,j+1}+r_{i-1,j}-r_{i-1,j+1} < p-n+j' \leq
r_{n,j+1}+r_{ij}-r_{i,j+1}$.  Using this it follows that the inverse
of $w_0^{(n)} \times \wh z(r)$ maps $\alpha^{-1}(N+n+2-p) = n+N+1 -
(p-n+j')$ to $n+N+1-z(r)(p-n+j') = n+N+1-(p-n+j'+\delta)$, and this
number furthermore lies in the interval $[n+r_{n,i+1}+1, n+r_{ni}]$.
Thus $\beta^{-1}$ subtracts $i'-1$, so the inverse of $\alpha
(w_0^{(n)} \times \wh z(r)) \beta$ maps $N+n+2-p$ to
$N+2n+2-p-\delta-i'-j'$, as required.
%
\end{proof}

We will finish this paper by proving that the double quiver polynomial
$K_r(a;b)$ is a specialization of the quiver formula of
\cite{buch:grothendieck}.  We need the following statement (about
power series in the variables $x_i = 1-a_i$ and $y_i = 1-b_i$.)

\begin{lemma} \label{L:super}
  Let $f(a;b) \in \Z[a_1^{-1},\dots,a_q^{-1}, b_1,\dots,b_p]$ be a
  Laurent polynomial that is separately symmetric in the variables
  $\{a_i\}$ and $\{b_i\}$.  Assume furthermore that
  $f(c,a_2,\dots,a_q; c,b_2,\dots,b_p)$ is independent of the variable
  $c$.  Then $f$ is a (possibly infinite) linear combination of double
  stable Grothendieck polynomials $G_\lambda(a;b)$.
\end{lemma}
\begin{proof}
  If we set $b_i = 1+y_i$ and $a_i = 1+x_i$, that is $(a_i)^{-1} =
  \sum_{k \geq 0} (-x_i)^k$, then the resulting power series $f(1+x_i;
  1+y_i)$ is {\em supersymmetric}, i.e.\ if one sets $x_1 = y_1$ then
  $f(1+x_i;1+y_i)$ becomes independent of these variables.  Since the
  lowest term of $G_\lambda(1+x_i; 1+y_i)$ is the double Schur
  polynomial $s_\lambda(x;y)$, it follows from
  \cite[Thm.~1]{stembridge:characterization} that there are
  coefficients $d_\lambda \in \Z$ such that the lowest term of
\[ f(1+x_i; 1+y_i) \, -\, \sum d_\lambda\, G_\lambda(1+x_i; 1+y_i) \]
has higher degree than the lowest term of $f(1+x_i;1+y_i)$.  The lemma
follows from this because each $G_\lambda(1+x_i;1+y_i)$ is also
supersymmetric \cite{fomin.kirillov:yang-baxter}.
\end{proof}

We will prove in \cite{buch:supersymmetry} that the linear combination
in this lemma is always finite, but this fact is not needed for the
proof of our last theorem.

\begin{thm} \label{T:double}
  For any rank conditions $r$ we have 
\[ K_r(a;b) = \sum_\mu c_\mu(r)\, G_{\mu_1}(a^1;b^0)\, 
   G_{\mu_2}(a^2;b^1) \cdots G_{\mu_n}(a^n;b^{n-1}) \,,
\]
where the sum is over all sequences $\mu = (\mu_1,\dots,\mu_n)$ of
partitions.
\end{thm}
\begin{proof}
  \refthm{T:multisuper} and \reflemma{L:super} imply that we can
  write
\[ K_{r+m}(a;b) \,=\, \sum c_\mu\, G_{\mu_1}(a^1;b^0)\, G_{\mu_2}(a^2;b^1)
  \cdots G_{\mu_n}(a^n;b^{n-1})
\]
where the sum is over (possibly infinitely many) sequences of
partitions $\mu$.  \refthm{T:rankstab} implies that the coefficients
$c_\mu$ are independent of $m$.  By setting $b^i = a^i$ for all $i$,
it therefore follows from \refthm{T:ratio} and the definition
(\ref{E:orig_kqf}) of quiver coefficients that $c_\mu = c_\mu(r)$ for
all $\mu$.  The result is obtained by setting $m=0$.
\end{proof}

By using \refcor{C:quivcoef}, the above theorem can also be
interpreted as a double component formula, that is
\[ K_r(a;b) = \sum_{(w_1,\dots,w_n)} (-1)^{\ell(w_1 \cdots w_n)-d(r)}\,
   G_{w_1}(a^1; b^0)\, G_{w_2}(a^2; b^1) \cdots G_{w_n}(a^n; b^{n-1}) \,,
\]
where the sum is over all KMS-factorizations for $r$.

\refthm{T:rankstab} makes it possible to extend the quiver polynomial
$K_r(a;b)$ to infinite sets of variables $a^i$ and $b^i$, by taking
the limit of $K_{r+m}(a;b)$ as $m$ tends to infinity.  Such limits are
called {\em double quiver functions\/} in
\cite{knutson.miller.ea:four}.  In the $K$-theory case it is
preferable to change to the variables $x^i_j = 1 - (a^i_j)^{-1}$ and
$y^i_j = 1 - b^i_j$ in order to obtain a nice formal power series.
This recovers the quiver formula $P_r$ constructed in \cite[\S
4]{buch:grothendieck}.  In fact, \refthm{T:double} shows that the
function $\lim_{m\to \infty} K_{r+m}(a;b)$ is obtained from $P_r$ by
setting $1^{\otimes i-1} \otimes x_j \otimes 1^{\otimes n-i} = 1 -
(a^i_j)^{-1}$ and $1^{\otimes i-1} \otimes y_j \otimes 1^{\otimes n-i}
= 1 - b^{i-1}_j$.  In particular, the cohomological double quiver
function used in \cite{knutson.miller.ea:four} is equal to the lowest
term of $P_r$.  Equivalently, this cohomological quiver function is a
specialization of the original quiver formula from
\cite{buch.fulton:chern} (see also the construction of this formula
given in \cite[\S 2]{buch:on}.)

By setting the variables $b^i_j$ equal to $1$ in \refthm{T:double},
one can deduce that general quiver coefficients are special cases of
the coefficients studied in \cite{buch.kresch.ea:grothendieck}.  More
details about this will be given in \cite{buch.sottile.ea:quiver},
together with a proof that quiver coefficients are special cases of
Schubert structure constants (see also
\cite{bergeron.sottile:schubert, lenart.robinson.ea:grothendieck}).




\providecommand{\bysame}{\leavevmode\hbox to3em{\hrulefill}\thinspace}
\providecommand{\MR}{\relax\ifhmode\unskip\space\fi MR }
\providecommand{\MRhref}[2]{%
  \href{http://www.ams.org/mathscinet-getitem?mr=#1}{#2}
}
\providecommand{\href}[2]{#2}

\end{document}